\documentclass[11pt]{dqaarticle}

\usepackage{
    amsmath,
    enumitem,
    hyperref,
    multicol,
    prftree,
    quiver,
    scalerel,
}

\usepackage[
    utf8,
]{inputenc}
\usepackage[
    T2A,
    T1,  
]{fontenc}
\usepackage{tgpagella}
\newcommand{\fallbackfont}[2]{%
    {\fontencoding{#1}\fontfamily{cmr}\selectfont#2}%
}
\definecolor{CiteGreen}{RGB}{0,150,0}
\hypersetup{
    colorlinks=true,
    citecolor=CiteGreen,
    filecolor=magenta,
    linkcolor=blue,
    urlcolor=cyan,
    pdfauthor={David Quinn ALVAREZ},
    pdftitle={A dual characterisation of simple and subdirectly-irreducible temporal Heyting algebras},
}
\newlist{axioms}{enumerate}{1}
\setlist[axioms]{
    align=left,
    labelwidth=\itemindent,
}

\def\runningtitle{A dual characterisation of simple and subdirectly-irreducible temporal Heyting algebras}
\title{%
    \runningtitle\\
    {\large \textit{with applications to the temporal Heyting calculus}}%
}
\def\runningauthor{D.\ Q.\ Alvarez}
\author{%
    David Quinn ALVAREZ\\
    {\small \texttt{dqalvarez@proton.me}}%
}
\date{10 June 2026}

\usepackage{
    adjustbox,
    amssymb,
    cancel,
    scalerel,
    stmaryrd,
}

\newcommand{\secref}[1]{\textsection\ref{#1}}

\newcommand{\notbin}[1]{\mathbin{\cancel{{#1}}}}
\newcommand{\axm}[1]{\textsf{(#1)}}
\newcommand{\axmref}[1]{\hyperref[#1]{\axm{#1}}}
\def\logdot{\ .\ }

\def\mimp{\Rightarrow}
\def\Mimp{\Longrightarrow}
\def\miff{\Leftrightarrow}
\def\Miff{\Longleftrightarrow}
\def\Diff{\mathbin{{:}{\Longleftrightarrow}}}
\def\Quad{\quad\quad\quad}
\def\LHS{\text{LHS}}
\def\RHS{\text{RHS}}

\def\contains{\ni}
\def\notcontains{\notbin{\ni}}
\def\empt{\varnothing}
\def\sub{\subseteq}
\def\sup{\supseteq}
\def\notsub{\notbin{\sub}}

\def\union{\cup}
\def\inter{\cap}
\def\Union{\bigcup}
\def\Inter{\bigcap}

\newcommand{\tup}[1]{{\langle #1 \rangle}}
\newcommand{\Tup}[1]{{\left\langle #1 \right\rangle}}

\def\fin{\textsf{fin}}
\def\pset{\mathop{{\wp}}}
\newcommand{\set}[1]{\{#1\}}
\newcommand{\Set}[1]{\left\{#1\right\}}

\newcommand{\bindset}[3]{\set{#1}_{#2}^{#3}}
\newcommand{\comp}[2]{\{ #1 \mid #2 \}}

\newcommand{\compin}[3]{\{ #1 \in #2 \mid #3 \}}

\def\too{\longrightarrow}
\def\mapstoo{\longmapsto}
\def\after{\mathbin{{\circ}}}
\def\ph{\raisebox{-1.5pt}{%
\begin{tikzpicture}
    \draw[densely dotted, line width=.6pt] (0,0) circle (4pt);
\end{tikzpicture}}
}

\def\id{\textsf{id}}
\newcommand{\idcat}[1]{\id_{\scaleobj{0.75}{#1}}}

\newcommand{\op}[1]{{#1^{\scaleobj{.8}{\textsf{op}}}}}
\def\iso{\cong}

\newcommand{\isocat}[1]{\mathbin{{\cong}^{\scaleobj{0.6}{#1}}}}
\def\equiv{\mathbin{{\simeq}}}

\def\ladj{\mathbin{{\dashv}}}

\def\Prop{\textsf{Prop}}
\DeclareRobustCommand\proves{\mathrel{|}\joinrel\mkern-.5mu\mathrel{-}}

\def\modelledby{\mathbin{\reflectbox{${\models}$}}}
\def\soco{\mathbin{{\modelledby}{\proves}}}
\def\phi{\varphi}
\def\imp{\rightarrow}
\def\biimp{\leftrightarrow}

\def\refutes{\notbin{\models}}

\def\Lm{\mathcal{L}_{\textsf{m}}}
\def\Lt{\mathcal{L}_{\textsf{t}}}

\let\oldBox\Box
\renewcommand{\Box}{\raisebox{-0.7pt}{${\oldBox}$}}
\def\bBox{\raisebox{-0.7pt}{${\blacksquare}$}}
\def\bDia{{\blacklozenge}}
\def\Dia{{\Diamond}}

\def\Term{\textsf{Term}}
\def\bDiaCom{\bDia\textsf{Com}}
\def\bDiaFilt{\bDia\Filt}
\def\Cong{\textsf{Cong}}
\def\Filt{\textsf{Filt}}
\def\fsi{\textsf{fsi}}
\def\Ideal{\textsf{Ideal}}
\def\PrFilt{\textsf{PrFilt}}

\newcommand{\filter}[1]{\left[\hspace{0.05cm}#1\right)}

\def\mydot{\scaleobj{0.75}{{\bullet}}}
\newcommand{\congfilt}[1]{#1_{\mydot}}
\newcommand{\filtcong}[1]{\mathbin{{#1^{\mydot}}}}

\def\Arc{\textsf{Arc}}
\def\ArcUp{\textsf{ArcUp}}
\def\ClArcUp{\textsf{ClArcUp}}
\def\ClDown{\textsf{ClDown}}

\def\ClopUp{\textsf{ClopUp}}
\def\ClUp{\textsf{ClUp}}

\def\OpUp{\textsf{OpUp}}
\def\ToRo{\textsf{ToRo}}
\def\tr{\mathbin{{\trianglelefteqslant}}}

\def\leq{\leqslant}
\def\geq{\geqslant}
\def\notleq{\notbin{\leq}}

\def\join{\vee}
\def\meet{\wedge}

\def\Meet{\bigwedge}

\newcommand{\upset}[1]{{\uparrow}#1}
\newcommand{\downset}[1]{{\downarrow}#1}
\def\Up{\textsf{Up}}

\def\Refl{\textsf{Refl}}

\def\Brel{\mathbin{B}}
\def\Rrel{\mathbin{R}}
\def\Zrel{\mathbin{Z}}
\def\thetarel{\mathbin{\theta}}
\def\Rb{\mathbin{R^{{\circ}}}}
\def\Rf{\mathbin{R}}

\def\Nats{\mathbb{N}}

\newcommand{\clop}[1]{#1^{*}}
\newcommand{\spec}[1]{#1_{*}}
\newcommand{\fram}[1]{#1_{{+}}}

\def\ES{\mathbf{ES}}
\def\fES{\mathbf{fES}}
\def\fHA{\mathbf{fHA}}
\def\HA{\mathbf{HA}}

\def\IPC{\mathbf{IPC}}
\def\mHC{\mathbf{mHC}}
\def\POS{\mathbf{POS}}
\def\tES{\mathbf{tES}}
\def\tHA{\mathbf{tHA}}
\def\tHC{\mathbf{tHC}}
\def\Tran{\mathbf{Tran}}
\def\tTran{\mathbf{tTran}}

\def\A{\mathbb{A}}
\def\B{\mathbb{B}}
\def\C{\mathbb{C}}
\def\K{\mathbf{K}}
\def\M{\mathbb{M}}
\def\V{\mathbf{V}}
\def\X{\mathbb{X}}
\def\Y{\mathbb{Y}}
\def\sepdot{\boldsymbol{{\cdot}}}
\def\R{\mathbin{R}}
\def\logdot{\ .\ }
\def\rbSig{\mathbin{r^{{\circ}}_{\Sigma}}}
\def\rfSig{\mathbin{r_{\Sigma}}}
\def\precSig{\mathbin{{{\prec}}_{\Sigma}}}
\def\RbSig{\mathbin{R^{{\circ}}_{\Sigma}}}
\def\RfSig{\mathbin{R_{\Sigma}}}
\def\leqSig{\mathbin{{{\leq}}_{\Sigma}}}
\def\nuSig{\nu_{\Sigma}}

\begin{document}

\maketitle
\thispagestyle{empty}

\begin{abstract}
    We establish an Esakia duality for the categories of temporal Heyting algebras and temporal Esakia spaces.
    This includes a proof of contravariant equivalence and a congruence/filter/closed-upset correspondence.
    We then study two notions of \mbox{« reachability »} on the relevant spaces/frames and show their equivalence in the finite case.
    We use these notions of reachability to give both lattice-theoretic and dual order-topological characterisations of simple and subdirectly-irreducible temporal Heyting algebras.
    Finally, we apply our duality results to prove the relational and algebraic finite model property for the temporal Heyting calculus.
    This, in conjunction with the proven characterisations, allows us to prove a relational completeness result that combines finiteness and the frame property dual to subdirect-irreducibility, giving us a class of finite, well-understood frames for the logic.
    ~
    \mbox{temporal Heyting calculus} $\sepdot$
    \mbox{Esakia duality} $\sepdot$
    \mbox{intuitionistic modal logic} $\sepdot$
    \mbox{temporal logic}
\end{abstract}

\section{Introduction}\label{sec:introduction}
    Temporal Heyting algebras were introduced in \cite{Esa06} as the algebraic models of the temporal Heyting calculus (also introduced in \cite{Esa06}).
    They contain, in addition to their Heyting algebra structure, the modal operators ${\bDia}$ and ${\Box}$, which, from a logical point of view, have the standard temporal reading.
    \begin{center}
    \renewcommand{\arraystretch}{1.33}
    \begin{tabular}{c|c|c}
        & Past & Future \\
        \hline
        ${\exists}$ & $\bDia$ & $\Dia$ \\
        \hline
        ${\forall}$ & $\bBox$ & $\Box$
    \end{tabular}
    \end{center}
    Note that, given the absence of negation on Heyting algebras, the modalities ${\bBox}$ and ${\Dia}$ are \textit{not} definable, making temporal Heyting algebras the most expressive models possible (given the ${\Box}$ of \cite{Esa06}) without entering into the \mbox{« jungle »} of ways of axiomatising all four modalities \cite{BouDieFerKre21}.
    Temporal Heyting algebras provide a natural setting to study the modality ${\bDia}$, the left adjoint of ${\Box}$, as algebraic adjointness has been studied extensively, especially in the context of modal logic \cite{MenSmi14}.
    The variety has been presented on at conferences \cite{Jib11,Als14}, but no proofs or papers have ever emerged.
    In particular, \cite{Als14} anticipated several of the results of the current paper, however, the central result therein contained a gap that requires correction (see \cite[Example 7.3.1]{Alv24}).
    \par
    The advantages of a duality-theoretic approach to modal and intuitionistic logic are well-known.
    Stone-like dualities allow us to approach logical questions from either algebraic or order-topological angles and give us the tools to transfer results between these two paradigms.
    In the former case, we have access to all of the tools of universal algebra and lattice theory.
    In the latter case, we can employ all of the familiar relational methods found in what is typically called \mbox{« modal} \mbox{logic »} as well as the methods of general topology.
    \par
    Of particular interest in algebraic logic are the simple and subdirectly-irreducible members of a given variety.
    From the perspective of universal algebra, these represent the most fundamental elements of the variety.
    In particular, the subdirectly-irreducible members will, in varietal settings, contain all of the truths of the corresponding logic, allowing us to restrict our study to only these members.
    For this reason, a thorough understanding of these members is essential to achieve a more general understanding of the variety and, in turn, the logic.
    \par
    In this paper, we develop an Esakia duality between the categories of algebraic and order-topological models of the temporal Heyting calculus.
    We give lattice-theoretic and dual order-topological characterisations of both simple and subdirectly-irreducible algebras.
    We then apply our duality results and characterisations to prove facts about the temporal Heyting calculus of more general interest.
    This includes the finite model property (FMP) and a relational completeness result that combines finiteness and the frame property dual to subdirect-irreducibility, giving us a class of finite, well-understood frames for the logic.
    \par
    Similar work has been done in a range of related settings.
    The most similar would be the work on \mbox{« Boolean} algebras with \mbox{operators »} (BAOs) in \cite{Ven04} and \mbox{« distributive} modal \mbox{algebras »} in \cite{Bir07}, though the results achieved in these two papers cannot be directly applied to our setting as our logic is sub-classical and the operators ${\rhd}$ and ${\lhd}$ of \cite{Bir07} are not definable on our algebras.
    Another similar paper is \cite{CasSagSan10}, in which an Esakia duality was developed for a reduct/supervariety of the current variety.
    Indeed, the current paper extends and expands upon the duality-theoretic results of \cite{CasSagSan10} by treating the temporal extension of the logic treated therein.
    Yet another similar paper is \cite{AkbAliMem24}+\cite{AkbAliMem24-2} (a two-part paper), which studies what is technically a generalisation of temporal Heyting algebras, but does not provide a dual characterisation of simple or subdirectly-irreducible algebras and gives a lattice-theoretic characterisation very different from the one contributed in the current paper.
    The current paper aims to demonstrate a combination of many of the results of the above-mentioned papers in a setting that has never been thoroughly treated in the literature.
    \par
    We note, for any interested reader, that we have authored an extensible symbolic model-checker for temporal Heyting algebras called \verb|thacheck|, which can be found at the following web address.
    \begin{center}
        \href{https://gitlab.com/dqalvarez/thacheck}{https://gitlab.com/dqalvarez/thacheck}
    \end{center}
    Will not discuss this software further in the current paper, but interested readers can refer to \cite[Appendix A]{Alv24} for an overview of the project.
    \par
    We provide the following section-by-section outline.
    \vspace{6pt}
    \begin{enumerate}[leftmargin=1.5cm,rightmargin=1.5cm]
        \item[\secref{sec:preliminaries}] We establish basic definitions of the relevant algebras, frames, and ordered topological spaces.
        We recall the basics of Esakia duality as well as an extension thereof.
        \item[\secref{sec:duality}] We define and establish an Esakia duality between the relevant categories of algebras and ordered topological spaces.
        We establish a congruence/filter/closed-upset correspondence.
        \item[\secref{sec:reachability}] We study two notions of \mbox{« reachability »} : one on ordered topological spaces and the other on finite underlying frames thereof.
        We show that these two notions are equivalent in the finite case.
        \item[\secref{sec:characterisations}] We give lattice-theoretic and order-topological characterisations of simple and subdirectly-irreducible algebras.
        In the finite case, we do the same element-wise and frame-theoretically.
        \item[\secref{sec:applications}] We apply our duality results to prove the relational and algebraic FMP.
        We use the algebraic FMP in conjunction with our characterisations to prove a final relational completeness result for the temporal Heyting calculus.
    \end{enumerate}
\section{Preliminaries}\label{sec:preliminaries}
    In this section, we define the relevant category/variety of algebras $\tHA$, class of frames (i.e.\ sets with binary relations) $\tTran$, and category of ordered topological spaces $\tES$ for the current paper.
    We do this by first recalling the definitions of intermediary categories/classes $\fHA$, $\Tran$, and $\fES$.
    (These are, themselves, modal augmentations of the classes of Heyting algebras \cite{Esa19,BezJon06}, posets \cite{BezJon06}, and Esakia spaces \cite{GehGoo24,Esa19} respectively.)
    We then recall some results of Esakia duality \cite{GehGoo24,Esa19} as well as a relevant specification to $\fHA$.
    \par
    First we establish some basic set-theoretic, algebraic, order-topological, and categorical notation.
    \begin{ntn}
        Let $X$ be a set, $w, w' \in X$, and $S,T \sub X$.
        We let $S - T$ denote set-theoretic difference, i.e.\ $\compin{x}{S}{x \notin T}$, and we let $-S$ denote set-theoretic complement, i.e.\ $X - S$.
        Given binary relations $R, R'$ on $X$, we define the following.
            \[ R^{-1} := \compin{\tup{y,x}}{X \times X}{\tup{x, y} \in R} \]
            \[ R[S] := \compin{y}{X}{\exists x \in S \logdot x \Rrel y} \Quad R[w] := R[\set{w}] \]
            \[ {R}{;}{R'} := \compin{\tup{x,z}}{X \times X}{\exists y \in X \logdot x \R y \R' z} \]
        Letting $\X := \tup{X, {R}}$, we denote the set of reflexive points on $X$ by $\Refl(\X)$.
        If ${\leq}$ is a partial order on $X$, then we let $\upset{S} := {\leq}[S]$ and $\downset{S} := {\leq}^{-1}[S]$ and $(w, w'] := \compin{x}{X}{w < x \leq w'}$.
        We also let the following denote the \mbox{« diagonal} \mbox{relation »} on $X$ (though we usually omit the subscript as it is clear from context) :
            \[ \Delta_{X} := \compin{\tup{x, y}}{X \times X}{ x = y} \]
        We let $\HA$ denote the category whose objects are Heyting algebras and whose morphisms are algebraic homomorphisms.
        Given $\A \in \HA$, we denote the set of \textit{congruences} of $\A$ by $\Cong(\A)$, the set of \textit{filters} of $\A$ by $\Filt(\A)$, the set of \textit{ideals} of $\A$ by $\Ideal(\A)$, and the set of \textit{prime filters} of $\A$ by $\PrFilt(\A)$.
        Given $S \sub \A$, we denote the filter generated by $S$ by
            \[ \filter{S} := \upset{\comp{a_1 \meet \dots \meet a_n}{\bindset{a_i}{i=1}{n} \sub S}}. \]
        We let $\ES$ denote the category whose objects are Esakia spaces and whose morphisms are Esakia morphisms (also called \mbox{« continuous} {p-morphisms »}).
        Given $\X \in \ES$, we denote the set of \textit{upsets} of $\X$ by $\Up(\X)$, the set of \textit{closed upsets} of $\X$ by $\ClUp(\X)$, the set of \textit{closed downsets} of $\X$ by $\ClDown(\X)$, the set of \textit{open upsets} of $\X$ by $\OpUp(\X)$, and the set of \textit{clopen upsets} of $\X$ by $\ClopUp(\X)$.
        Given a class of structures $\K$, we denote the class of finite members of $\K$ by $\K_{\fin}$.
        If $\K$ is a class of algebras, we denote the class of finite subdirectly-irreducible members of $\K$ by $\K_{\fsi}$.
        Finally, given a category $\mathbf{C}$ and a $\mathbf{C}$-object $\X$, we denote the identity morphism on $\X$ by $\idcat{\X}$.
    \end{ntn}
    \par
    We now recall the definition of \textit{frontal Heyting algebras}, defined in \cite{Esa06} and studied extensively in \cite{CasSagSan10}.
    \begin{dfn}[Frontal Heyting algebra]
        A \textit{frontal Heyting algebra} is an algebra $\tup{A, {\meet}, {\join}, {\imp}, {\Box}, 0, 1}$ such that $\tup{A, {\meet}, {\join}, {\imp}, 0, 1}$ is a Heyting algebra and for all $a,b \in A$,
            \[ \Box(a \meet b) = \Box a \meet \Box b \Quad a \leq \Box a \Quad \Box a \leq b \join (b \imp a). \]
        The category $\fHA$ has as its objects frontal Heyting algebras and as its morphisms algebraic homomorphisms.
    \end{dfn}
    We now define the main variety/category of the current paper.
    \begin{dfn}[Temporal Heyting algebra]
        A \textit{temporal Heyting algebra} is an algebra $\tup{A, {\meet}, {\join}, {\imp}, {\bDia}, {\Box}, 0, 1}$ such that $\tup{A, {\meet}, {\join}, {\imp}, {\Box}, 0, 1}$ is a frontal Heyting algebra and ${\bDia} \ladj {\Box}$, i.e.
            \[ \bDia a \leq b \quad \Miff \quad a \leq \Box b. \]
        The category $\tHA$ has as its objects temporal Heyting algebras and as its morphisms algebraic homomorphisms.
    \end{dfn}
    Note that the following (in)equalities yield an equational definition of temporal Heyting algebras :
        \[ \bDia 0 = 0 \Quad \bDia (a \join b) = \bDia a \join \bDia b \Quad a \leq \Box\bDia a \Quad \bDia \Box a \leq a. \]
    This is mentioned to ensure the reader that $\tHA$ forms a \textit{variety}, which will be relevant when establishing Theorem \ref{thm:thc-soco-tha-fsi} much later on.
    \par
    Having defined our algebras, we define the underlying frames of the dual ordered topological spaces.
    \begin{dfn}[Transit]
        A \textit{transit} ($\Tran$) is a frame $\tup{X, R, {\leq}}$ such that $\tup{X, {\leq}}$ is a poset and ${\leq}$ is the reflexivisation of $R$, i.e.\ ${\leq} = R \union \Delta_{X}$ or, equivalently,
            \[ x \leq y \quad \Miff \quad x \R y \text{ or } x = y. \]
    \end{dfn}
    Note that it can easily be shown that the relation $R$ on transits is antisymmetric, transitive, and satisfies the \mbox{« mix »} condition $R = {\leq}{;}{R}{;}{\leq}$ as well as the inequalities ${<} \sub {R}$ and ${R} \sub {\leq}$.
    This relation can be thought of as identical to the ${\leq}$ relation except for the fact that it could be missing some (or, potentially, all) reflexive loops.
    \begin{exm}\label{exm:transit}
        The following is an example of a transit.
        
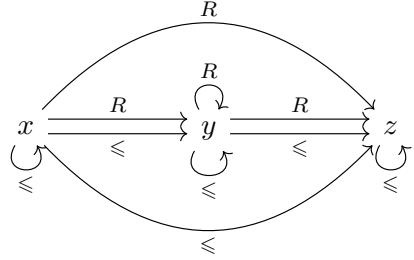
\[\begin{tikzcd}[ampersand replacement=\&]
            x \&\& y \&\& z
            \arrow["{{\leq}}"', from=1-1, to=1-1, loop, in=300, out=240, distance=5mm]
            \arrow["R", shift left, from=1-1, to=1-3]
            \arrow["{{\leq}}"', shift right, from=1-1, to=1-3]
            \arrow["R", curve={height=-50pt}, from=1-1, to=1-5]
            \arrow["{{\leq}}"', curve={height=50pt}, from=1-1, to=1-5]
            \arrow["R", from=1-3, to=1-3, loop, in=60, out=120, distance=5mm]
            \arrow["{{\leq}}"', from=1-3, to=1-3, loop, in=300, out=240, distance=5mm]
            \arrow["R", shift left, from=1-3, to=1-5]
            \arrow["{{\leq}}"', shift right, from=1-3, to=1-5]
            \arrow["{{\leq}}"', from=1-5, to=1-5, loop, in=300, out=240, distance=5mm]
        \end{tikzcd}\]
        Observe that the relations $R$ and ${\leq}$ agree everywhere except for the reflexive loops : $R$ is missing $\tup{x, x}$ and $\tup{z, z}$.
        This, however, implies that the reflexivisation of $R$ is identical to ${\leq}$ as desired.
        For this reason, it suffices to depict only the $R$ relation, as the ${\leq}$ relation can be inferred.
        Furthermore, the transitive arrows are typically omitted, leaving us with the following depiction of $R$ from which the above diagram can be inferred.
        
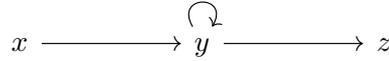
\[\begin{tikzcd}[ampersand replacement=\&]
            x \&\& y \&\& z
            \arrow[from=1-1, to=1-3]
            \arrow[from=1-3, to=1-3, loop, in=60, out=120, distance=5mm]
            \arrow[from=1-3, to=1-5]
        \end{tikzcd}\]
    \end{exm}
    \begin{dfn}[Temporal transit]
        A \textit{temporal transit} ($\tTran$) is a frame $\tup{X, {\Rb}, {\Rf}, {\leq}}$ such that $\tup{X, {\Rf}, {\leq}}$ is a transit and ${\Rb}$ is the inverse of ${\Rf}$, i.e.\ $x \Rf y$ iff $y \Rb x$.
    \end{dfn}
    \begin{rmk}
        The reader may have noticed that temporal transits could be equivalently defined using only a « forward-facing » relation $R$.
        The poset relation ${\leq}$ could then be defined as the reflexivisation of $R$ and the « backward-facing » relation could be defined as $R^{-1}$.
        Indeed, this was the way transits were defined and treated in \cite{Esa06,CasSagSan10}.
        However, we have decided that as the dual operations ${\imp}$, ${\Box}$, ${\bDia}$ are all explicit in the signature of the dual algebras and in the language of the relevant logic (to-be-defined in \secref{sec:applications}), the proofs read more easily when all three relations are made explicit and included in the signature.
    \end{rmk}
    \par
    We now define the Esakia spaces with underlying (temporal) transits that will be dual to the above-defined algebras.
    \begin{dfn}[Frontal Esakia space]
        A \textit{frontal Esakia space} (called an \mbox{« $Rf$}-Heyting \mbox{space »} or, simply, an \mbox{« $Rf$-space »} in \cite{CasSagSan10}) is an ordered topological space $\X := \tup{X, R, {\leq}, \Omega}$ such that $\tup{X, {\leq}, \Omega}$ is an Esakia space and for all $x \in X$ and $K \sub X$,
        \begin{itemize}
            \item $\tup{X, {\Rrel}, {\leq}}$ is a transit
            \item $K \in \ClopUp(\X)$ implies $-{R^{-1}}[-K] \in \ClopUp(\X)$
            \item $R[x]$ is closed.
        \end{itemize}
        A \textit{frontal Esakia morphism} (called an \mbox{« $Rf$-morphism »} in \cite{CasSagSan10}) is a map $f : \X \to \Y$ such that $f$ is an Esakia morphism and for all $x_1, x_2 \in \X$,
        \begin{itemize}
            \item $x_1 \Rrel x_2$ implies $f x_1 \Rrel f x_2 $
            \item $f x_1 \Rrel y$ implies $\exists x_3 \in \X$ such that $x_1 \Rrel x_3$ and $f x_3 = y$.
        \end{itemize}
        The category $\fES$ has as its objects frontal Esakia spaces and as its morphisms frontal Esakia morphisms.
    \end{dfn}
    \begin{rmk}
    The constraint that $R[x]$ be closed was omitted in \cite{CasSagSan10}, but is necessary in analogous contexts to prove that the natural isomorphism in the topological category $\gamma_{\X} : \X \to \spec{\clop{\X}}$ (to-be-defined) has the back-and-forth p-morphism conditions.
    The constraint is easily proven to be true on the dual spaces of $\fHA$s, but does not follow from the other constraints, so it was likely a small oversight, bearing no meaningful consequences to the duality results of \cite{CasSagSan10}.
    \end{rmk}
    \begin{dfn}[Temporal Esakia space]
        A \textit{temporal Esakia space} is an ordered topological space $\X := \tup{X, {\Rb}, {\Rf}, {\leq}, \Omega}$ such that $\tup{X, {\Rf}, {\leq}, \Omega}$ is a frontal Esakia space and for all $x \in X$ and $K \sub X$,
        \begin{axioms}
            \item[\axm{tES.o.1}]\label{tES.o.1} $\tup{X, {\Rb}, {\Rf}, {\leq}}$ is a temporal transit
            \item[\axm{tES.o.2}]\label{tES.o.2} $K \in \ClopUp(\X)$ implies $({\Rb})^{-1}[K] \in \ClopUp(\X)$
            \item[\axm{tES.o.3}]\label{tES.o.3} ${\Rb}[x]$ is closed.
        \end{axioms}
        A \textit{temporal Esakia morphism} is a map $f : \X \to \Y$ such that $f$ is a frontal Esakia morphism and for all $x_3, x_2 \in \X$ and $y \in \Y$,
        \begin{axioms}
            \item[\axm{tES.m.1}]\label{tES.m.1} $x_3 \Rb x_2$ implies $f x_3 \Rb f x_2$.
            \item[\axm{tES.m.2}]\label{tES.m.2} $f x_2 \Rb y$ implies $\exists x_1 \in \X$ such that $x_2 \Rb x_1$ and $y \leq f x_1$.
        \end{axioms}
        The category $\tES$ has as its objects temporal Esakia spaces and as its morphisms temporal Esakia morphisms.
    \end{dfn}
    Note that the constraint \axmref{tES.m.1} follows from the morphism conditions of $\fES$ combined with the definition of temporal transits.
    Note, also, that in \axmref{tES.o.2}, $({\Rb})^{-1}[K]$ is equal to ${\Rf}[K]$.
    Though these choices sacrifice minimality and some simplicity of the definition, the constraints were so-included for explicit agreement with typical axiomatisations of p-morphisms.
    \par
    The following table is included to help the reader navigate these definitions.
    The final row, in particular the variety $\tHA$, is the central focus of the current paper.
    Note that the \mbox{« Language »} column lists connectives \textit{in addition to} ${\meet}, {\join}, \bot, \top$.
    \vspace{4pt}
    \begin{center}
    \renewcommand{\arraystretch}{1.33}
    \begin{tabular}{ |c|c|c|c|c|c| }
    \hline
    Language & References & Logic & Algebras & Frames & Spaces \\
    \hline
    \hline
    ${\imp}$ & \cite{BezJon06, ChaZak97,Esa19} & $\IPC$ & $\HA$ & $\POS$ & $\ES$ \\
    & & \rotatebox[origin=c]{270}{$\sub$} & \rotatebox[origin=c]{90}{$\sub$} & \rotatebox[origin=c]{90}{$\sub$} & \rotatebox[origin=c]{90}{$\sub$} \\
    ${\imp}, {\Box}$ & \cite{Esa06,CasSagSan10} & $\mHC$ & $\fHA$ & $\Tran$ & $\fES$ \\
    & & \rotatebox[origin=c]{270}{$\sub$} & \rotatebox[origin=c]{90}{$\sub$} & \rotatebox[origin=c]{90}{$\sub$} & \rotatebox[origin=c]{90}{$\sub$} \\
    ${\imp}, {\Box}, {\bDia}$ & \cite{Esa06,Alv24} & $\tHC$ & $\boxed{\tHA}$ & $\tTran$ & $\tES$ \\
    \hline
    \end{tabular}
    \end{center}
    \vspace{6pt}
    \par
    We now recall the basics of Esakia duality.
    To do so we first establish some notation.
    Given $\A \in \HA$, recall the map
        \[ \pi_{\A} : \A \too \pset(\PrFilt(\A)) \Quad \pi_{\A} : a \mapstoo \compin{x}{\PrFilt(\A)}{a \in x}.\]
    We let $\spec{\ph} : \HA \leftrightarrows \ES : \clop{\ph}$ be the pseudo-inverse functors witnessing the contravariant equivalence $\HA \simeq \op{\ES}$, known as \mbox{« Esakia} \mbox{duality »}, established in \cite{Esa74} and treated extensively in \cite{Esa19}.
    Recall that
        \[ \spec{\A} := \tup{\PrFilt(\A), {\sub}, \Omega_{\A}} \]
    where $\Omega_{\A}$ is the topology generated by the subbasis
        \[ \comp{\pi a}{a \in \A} \union \comp{-\pi a}{a \in \A} \]
    and that $\spec{h} := h^{-1}[\ph]$.
    Recall, also, that
        \[ \clop{\X} := \tup{\ClopUp(\X), {\inter}, {\union}, {\imp}, \empt, \X} \]
    where
        \[ {K_1 \imp K_2} := -\downset{(K_1 - K_2)} \]
    and that $\clop{f} := f^{-1}[\ph]$.
    Furthermore, we let $\pi : \idcat{\HA} \Rightarrow \clop{\spec{\ph}}$ and $\gamma : \idcat{\ES} \Rightarrow \spec{\clop{\ph}}$ be the natural isomorphisms witnessing the equivalence in $\HA$ and $\ES$ respectively, implying that the following squares commute for all $\HA$-morphisms $h : \A \to \B$ and all $\ES$-morphisms $f: \X \to \Y$.
    \[\begin{tikzcd}[ampersand replacement=\&]
        \A \&\& {\clop{\spec{\A}}} \&\& \X \&\& {\spec{\clop{\X}}} \\
        \B \&\& {\clop{\spec{\B}}} \&\& \Y \&\& {\spec{\clop{\Y}}}
        \arrow["{\pi_{\A}}", tail, two heads, from=1-1, to=1-3]
        \arrow["h"', from=1-1, to=2-1]
        \arrow["{\clop{\spec{h}}}", from=1-3, to=2-3]
        \arrow["{\gamma_{\X}}", tail, two heads, from=1-5, to=1-7]
        \arrow["f"', from=1-5, to=2-5]
        \arrow["{\spec{\clop{f}}}", from=1-7, to=2-7]
        \arrow["{\pi_{\B}}"', tail, two heads, from=2-1, to=2-3]
        \arrow["{\gamma_{\Y}}"', tail, two heads, from=2-5, to=2-7]
    \end{tikzcd}\]
    Recall that $\gamma_{\X}$ is defined by the rule $x \mapsto \compin{K}{\ClopUp(\X)}{x \in K}$.
    When clear from context, we omit the subscripts of $\pi$ and $\gamma$.
    We also rely on the reader to parse $\clop{\spec{\ph}}$ as only one of $\clop{(\spec{\ph})}$ and $\spec{(\clop{\ph})}$ will be well-defined for a given operand.
    \par
    Recall also the class operator $\fram{\ph} : \HA \to \POS$ which maps an algebra $\A$ to the underlying frame of $\spec{\A}$ (forgetting the topology).
    \par
    We also recall the following poset-isomorphisms for $\A \in \HA$ between congruences, filters, and closed upsets respectively.
        \[ \tup{\Cong(\A), {\sub}} \isocat{\POS} \tup{\Filt(\A), {\sub}} \isocat{\POS} \tup{\ClUp(\spec{\A}), {\sup}} \]
    The maps are as follows where $\congfilt{\theta} := [1]_{\theta}$ and ${\filtcong{F}} := \compin{\tup{a, b}}{\A \times \A}{a \biimp b \in F}$.
    \[\begin{tikzcd}[ampersand replacement=\&]
        {\tup{\Cong(\A), {\sub}}} \&\& {\tup{\Filt(\A), {\sub}}} \&\& {\tup{\ClUp(\spec{\A}), {\sup}}}
        \arrow["{\congfilt{\ph}}", shift left=2, tail, two heads, from=1-1, to=1-3]
        \arrow["{{\filtcong{\ph}}}", shift left=2, tail, two heads, from=1-3, to=1-1]
        \arrow["{\Inter\pi[\ph]}", shift left=2, tail, two heads, from=1-3, to=1-5]
        \arrow["\Inter", shift left=2, tail, two heads, from=1-5, to=1-3]
    \end{tikzcd}\]
    Recall also that $\Union\pi[\ph]$ is a poset-isomorphism $\tup{\Ideal(\A), {\sub}} \to \tup{\OpUp(\A), {\sub}}$ and that for all $F \in \Filt(\A)$ and $I \in \Ideal(\A)$,
        \[ x \in \Inter\pi[F] \quad \Miff \quad F \sub x \Quad x \in \Union\pi[I] \quad\Miff\quad x \inter I \neq \empt. \]
    \par
    The functors $\spec{\ph} : \HA \leftrightarrows \ES : \clop{\ph}$ were modified by \cite{CasSagSan10} in the following way to be defined $\fHA \leftrightarrows \fES$.
    Given $\A \in \fHA$, the relation $R$ was defined on the dual space $\spec{\A}$ as
        \[ x \Rrel y \quad \Diff \quad \forall a \in \A \logdot \Box a \in x \mimp a \in y. \]
    Given $\X \in \fES$, the operation ${\Box}$ was defined on the dual algebra $\clop{\X}$ as
        \[ \Box K := -{R^{-1}}[-K] = \compin{x}{\X}{R[x] \sub K}. \]
    Given these modifications, it was shown in \cite{CasSagSan10} that $\fHA \equiv \op{\fES}$.
\section{Esakia duality for tHA}\label{sec:duality}
    In this section, we establish an Esakia duality between the categories $\tHA$ and $\tES$.
    We then establish a congruence/filter/closed-upset correspondence that will be used extensively throughout the subsequent sections.
    \par
    We begin by further modifying our functors $\spec{\ph} : \fHA \leftrightarrows \fES : \clop{\ph}$ so that they are defined $\tHA \leftrightarrows \tES$.
    Given $\A \in \tHA$, the relation ${\Rb}$ is defined on the dual space $\spec{\A}$ as
        \[ x \Rb w \quad \Diff \quad \forall a \in \A \logdot a \in w \mimp \bDia a \in x. \]
    (Note that we now have $\fram{\ph} : \tHA \to \tTran$.)
    Given $\X \in \tES$, the operation ${\bDia}$ is defined on the dual algebra $\clop{\X}$ as
        \[ \bDia K := ({\Rb})^{-1}[K] = {\Rf}[K] = \compin{x}{\X}{\Rb[x] \inter K \neq \empt}. \]
    \par
    We show the well-definedness of $\spec{\ph} : \tHA \leftrightarrows \tES : \clop{\ph}$ for objects and morphisms in the following lemmas (Lemmas \ref{lem:spec-well-defined-tha-obj} through \ref{lem:clop-well-defined-tes-mor}), but first we recall the \textit{Prime Filter Theorem}.
    \begin{fct}[Prime Filter Theorem]
        Given $\A \in \HA$, a filter $F \in \Filt(\A)$, and an ideal $I \in \Ideal(\A)$ such that $F \inter I = \empt$, there exists a prime filter $x \in \PrFilt(\A)$ such that $F \sub x$ and $x \inter I = \empt$.
    \end{fct}
    \noproof
    \begin{lem}\label{lem:spec-well-defined-tha-obj}
        Given $\A \in \tHA$, we have $\spec{\A} \in \tES$.
    \end{lem}
    \begin{prf}
        Let $\X := \spec{\A}$.
        ~\\~\\
        \axmref{tES.o.1}
            Assuming $x \Rf y$, we show $y \Rb x$.
            Given $a \in x$, we have $\bDia a \leq \bDia a$, implying $a \leq \Box \bDia a$, further implying that $\Box \bDia a \in x$.
            Since $x \Rf y$, we have $\bDia a \in y$.
            Now, assuming $y \Rb x$, we show $x \Rf y$.
            Given $a \notin y$, we have $\Box a \leq \Box a$, implying $\bDia \Box a \leq a$, further implying that $\bDia \Box a \notin y$.
            Since $y \Rb x$, we have $\Box a \notin x$.
        ~\\~\\
        \axmref{tES.o.2}
            Given $K \in \ClopUp(\X)$, we know, since $\pi_{\A}$ is surjective, that there is some $a \in \A$ such that $K = \pi a$.
            We claim that ${\Rf}[K] = \pi \bDia a$.
            Note that if this is the case, then ${\Rf}[K] \in \ClopUp(\X)$ as desired.
            $(\sub)$
                Given $x \in {\Rf}[K]$, we have $w \in K$ and $w \Rf x$, implying that $x \Rb w$.
                Since $w \in \pi a$, we have $a \in w$, implying $\bDia a \in x$, which implies that $x \in \pi \bDia a$.
            $(\sup)$
                Given $x \in \pi \bDia a$, we find a prime filter $w$ such that $w \in \pi a = K$ and $w \Rf x$.
                Consider the filter $\upset{a}$ and the ideal $\bDia^{-1}[-x]$.
                (It can easily be shown that $\bDia^{-1}[-x]$ contains $0$ and is, therefore, non-empty, is downwards-closed, and is ${\join}$-closed, implying that it is, indeed, an ideal.)
                Assume, toward a contradiction, that these sets are \textit{not} disjoint.
                Then there is some $b \in \upset{a} \inter \bDia^{-1}[-x]$, implying that $a \leq b$ and $\bDia b \in -x$.
                Since $a \leq b$, we have $\bDia a \leq \bDia b$, so, since $-x$ is a downset, we have $\bDia a \in -x$, implying $\bDia a \notin x$, contradicting our assumption that $x \in \pi \bDia a$.
                So, since $\upset{a} \inter \bDia^{-1}[-x] = \empt$, we know, by the Prime Filter Theorem, that there exists some $w \in \PrFilt(\A)$ such that $\upset{a} \sub w$ and ${w \inter \bDia^{-1}[-x] = \empt}$.
                Since $\upset{a} \sub w$, we have $w \in \pi a = K$.
                To see that $x \Rb w$, observe that given $\bDia b \notin x$, we have $b \in \bDia^{-1}[-x]$, implying $b \notin w$.
                So $w \Rf x$ and $w \in K$, implying $x \in {\Rf}[K]$.
        ~\\~\\
        \axmref{tES.o.3}
            Given $x \in \X$, we claim that $\Rb[x] = -\Union\pi[\bDia^{-1}[-x]]$.
            Note that if this is the case, then since $\bDia^{-1}[-x] \in \Ideal(\A)$, we will have $\Union\pi[\bDia^{-1}[-x]] \in \OpUp(\X)$, implying that $\Rb[x] = -\Union\pi[\bDia^{-1}[-x]] \in \ClDown(\X)$.
            ~\\
            \begin{align*}
                w \in \Rb[x] \quad &\Miff \quad x \Rb w \\
                &\Miff \quad \forall a \in \A \logdot a \in w \mimp \bDia a \in x \\
                &\Miff \quad \neg \exists a \in \A \logdot a \in w \text{ and } \bDia a \in - x \\
                &\Miff \quad \neg \exists a \in \A \logdot a \in w \text{ and } a \in \bDia^{-1}[- x] \\
                &\Miff \quad w \inter \bDia^{-1}[- x] = \empt \\
                &\Miff \quad w \notin \Union\pi[\bDia^{-1}[-x]] \\
                &\Miff \quad w \in - \Union\pi[\bDia^{-1}[-x]]
            \end{align*}
    \end{prf}
    \begin{lem}\label{lem:spec-well-defined-tha-mor}
        Given a $\tHA$-morphism $h : \A \to \B$, the map $\spec{h} : \spec{\B} \to \spec{\A}$ is a $\tES$-morphism.
    \end{lem}
    \begin{prf}
        Let
            \[ \X := \spec{\B} \Quad \Y := \spec{\A} \Quad f :=  \spec{h} : \X \to \Y. \]
        \axmref{tES.m.1}
            Given $x_3, x_2 \in \X$ such that $x_3 \Rb x_2$, we show that $f x_3 \Rb f x_2$.
            But this is quite obvious as $x_3 \Rb x_2$ implies $x_2 \Rf x_3$, implying, since $f$ is an $\fES$-morphism, that $f x_2 \Rf f x_3$, finally implying that $f x_3 \Rb f x_2$.
        ~\\~\\
        \axmref{tES.m.2}
            Given $y \in \Y$ and $x_2 \in \X$ such that $f x_2 \Rb y$, we show that there must exist some $x_1 \in \X$ such that $x_2 \Rb x_1$ and $y \sub f x_1$.
            \[\begin{tikzcd}[ampersand replacement=\&]
                {x_2} \&\&\& {f x_2} \\
                {x_1} \&\& {f x_1} \\
                \&\&\& y
                \arrow["f", from=1-1, to=1-4]
                \arrow["{{\Rb}}"', dashed, from=1-1, to=2-1]
                \arrow["{{\Rb}}", from=1-4, to=3-4]
                \arrow["f"', from=2-1, to=2-3]
                \arrow["{{\sub}}", dashed, from=3-4, to=2-3]
            \end{tikzcd}\]
            Consider the filter $\upset{h[y]}$ and the ideal $\bDia^{-1}[-x_2]$.
            Assume, toward a contradiction, that these sets are \textit{not} disjoint.
            Then there is some $b \in \upset{h[y]} \inter \bDia^{-1}[-x_2]$, implying that there exists some $a \in y$ such that $ha \leq b$ and $\bDia b \notin x_2$.
            Now since $ ha \leq b$, we have $h \bDia a = \bDia ha \leq \bDia b$ (using that fact that $h$ is homomorphic over ${\bDia}$).
            Also, since $a \in y$ and $fx_2 \Rb y$, we have $\bDia a \in f x_2 = h^{-1}[x_2]$, implying that $h\bDia a \in x_2$.
            But this implies that $\bDia b \in x_2$, contradicting our assumption that $\bDia b \notin x_2$.
            So, since $\upset{h[y]} \inter \bDia^{-1}[-x_2] = \empt$, we know, by the Prime Filter Theorem, that there exists some $x_1 \in \spec{\B}$ such that $\upset{h[y]} \sub x_1$ and $x_1 \inter \bDia^{-1}[-x_2] = \empt$.
            Now we claim that $x_2 \Rb x_1$ and $y \sub fx_1$.
            To see that $x_2 \Rb x_1$, consider some $\bDia b \notin x_2$ and observe that $b \in \bDia^{-1}[-x_2]$, implying that $b \notin x_1$.
            To see the $y \sub fx_1$, consider some $a \in y$ and observe that $ha \in h[y] \sub \upset{h[y]} \sub x_1$, implying $a \in h^{-1}[x_1] = f x_1$.
    \end{prf}
    \begin{lem}\label{lem:clop-well-defined-tes-obj}
        Given $\X \in \tES$, we have $\clop{\X} \in \tHA$.
    \end{lem}
    \begin{prf}
        Let $\A := \clop{\X}$.
        To address the well-definedness of ${\bDia}$ on $\A$, note that $K \in \A$ implies $\bDia{K} = {\Rf}[K] \in \A$, implying that $\bDia: \A \to \A$.
        It remains to show that for $K_1, K_2 \in \A$, we have
            \[ \bDia K_1 \sub K_2 \quad \Miff \quad K_1 \sub \Box K_2. \]
        $(\Rightarrow)$
            We let $x \in K_1$ and we show that $x \in \Box K_2$, i.e.\ $\Rf[x] \sub K_2$.
            Given $y \in \Rf[x]$, we have $x \Rf y$, implying $y \Rb x$, further implying, since $x \in K_1$, that $\Rb[y] \inter K_1 \neq \empt$, finally implying that $y \in \bDia K_1$.
            By the $\LHS$, this implies that $y \in K_2$ as desired.
        ~\\~\\
        $(\Leftarrow)$
            We let $x \in \bDia K_1$, i.e.\ $\Rb[x] \inter K_1 \neq \empt$, and show that $x \in K_2$.
            Given $\Rb[x] \inter K_1 \neq \empt$, we have $x \Rb w$ and $w \in K_1$, implying, by the $\RHS$, that $w \in \Box K_2$, finally implying that $\Rf[w] \sub K_2$.
            Since $x \in \Rf[w]$, we have $x \in K_2$ as desired.
    \end{prf}
    \begin{lem}\label{lem:clop-well-defined-tes-mor}
        Given a $\tES$-morphism $f : \X \to \Y$, the map $\clop{f} : \clop{\Y} \to \clop{\X}$ is a $\tHA$-morphism.
    \end{lem}
    \begin{prf}
        Let
            \[ \A := \clop{\Y} \Quad \B := \clop{\X} \Quad h := \clop{f} : \A \to \B.\]
        We show that $h$ is homomorphic over ${\bDia}$, i.e., given $K \in \A$, we have
            \[ h \bDia K = \bDia h K. \]
        $(\sub)$
            Let $x_2 \in h\bDia K = f^{-1}[\bDia K]$, implying that $f x_2 \in \bDia K$.
            This implies there is some $y \in \Y$ such that $y \in {\Rb}[f x_2] \inter K$, implying that $f x_2 \Rb y$.
            By \axmref{tES.m.2}, there exists some $x_1 \in \X$ such that $x_2 \Rb x_1$ and $y \leq f x_1$.
            Also, since $y \leq f x_1$ and $y \in K$, we have $f x_1 \in K$, implying $x_1 \in f^{-1}[K] = hK$.
            So since $x_1 \in {\Rb}[x_2] \inter hK$, we have $x_2 \in \bDia h K$.
        ~\\~\\
        $(\sup)$
            Given $x_2 \in \bDia h K$, we have ${\Rb}[x_2] \inter hK \neq \empt$, implying $x_2 \Rb x_1$ and $x_1 \in hK = f^{-1}[K]$, implying $f x_1 \in K$.
            Since $x_2 \Rb x_1$, \axmref{tES.m.1} implies that $f x_2 \Rb f x_1$, implying $f x_1 \in {\Rb}[f x_2] \inter K$.
            This implies that $f x_2 \in \bDia K$, finally implying that $x_2 \in f^{-1}[\bDia K] = h \bDia K$.
    \end{prf}
    \par
    We now show that $\pi$ and $\gamma$ are well-defined natural isomorphisms.
    \begin{lem}\label{lem:pi-well-defined-nat-trans}
        Given $\A \in \tHA$, the map $\pi_{\A} : \A \to \clop{\spec{\A}}$ is a $\tHA$-isomorphism.
    \end{lem}
    \begin{prf}
        Since it is known from the duality established in \cite{CasSagSan10} that $\pi$ is an $\fHA$-isomorphism, it suffices to show that $\pi$ is homomorphic over ${\bDia}$, i.e., given $a \in \A$, we have $\pi \bDia a = \bDia \pi a$.
        ~\\~\\
        $(\sup)$
            Given $x \in \bDia \pi a$, we have $\Rb[x] \inter \pi a \neq \empt$, so we have $x \Rb w$ and $w \in \pi a$.
            This implies that $a \in w$ implying that $\bDia a \in x$, finally implying that $x \in \pi \bDia a$.
        ~\\~\\
        $(\sub)$
            Given $x \in \pi \bDia a$, we have $\bDia a \in x$.
            It can be easily checked that $\bDia^{-1}[-x] \in \Ideal(\A)$.
            Now we claim that $\upset{a} \inter \bDia^{-1}[-x] = \empt$.
            This must be the case as if it were not, then we would have $a \leq b$ and $\bDia b \notin x$, implying that $\bDia a \leq \bDia b$, further implying that $\bDia b \in x$, giving us a contradiction.
            Since $\upset{a} \inter \bDia^{-1}[-x] = \empt$, the Prime Filter Theorem implies that there exists some $w \in \spec{\A}$ such that $\upset{a} \sub w$ and $w \inter \bDia^{-1}[-x]= \empt$.
            Now we claim that $w \in \Rb[x] \inter \pi a$.
            Since $\upset{a} \sub w$, we have $a \in w$, implying $w \in \pi a$.
            And given $b \in w$, we have $b \notin \bDia^{-1}[-x]$, implying $\bDia b \notin -x$, implying $\bDia b \in x$.
            This shows that $x \Rb w$.
            So since $\Rb[x] \inter \pi a \neq \empt$, we have $x \in \bDia \pi a$ as desired.
    \end{prf}
    \begin{lem}\label{lem:gamma-well-defined-nat-trans}
        Given $\X \in \tES$, the map $\gamma_{\X} : \X \to \spec{\clop{\X}}$ is a $\tES$-isomorphism.
    \end{lem}
    \begin{prf}
        Since it is known from the duality established in \cite{CasSagSan10} that $\gamma$ is an $\fES$-isomorphism, it suffices to show that $\gamma$ is a $\tES$-morphism.
        ~\\~\\
        \axmref{tES.m.1}
            Given $x_3, x_2 \in \X$ such that $x_3 \Rb x_2$, we show $\gamma x_3 \Rb \gamma x_2$.
            Given $K \in \gamma x_2$, we have $x_2 \in K$, implying that $x_2 \in \Rb[x_3] \inter K$.
            This implies that $x_3 \in \bDia K$, implying $\bDia K \in \gamma x_3$ as desired.
        ~\\~\\
        \axmref{tES.m.2}
            Given $x_2 \in \X$ and $y \in \spec{\clop{\X}}$ such that $\gamma x_2 \Rb y$, and show that there must exist some $x_1 \in \X$ such that $x_2 \Rb x_1$ and $y \sub \gamma x_1$.
            Assume, toward a contradiction, that no such $x_1$ exists, i.e.
                \[ \neg \exists x_1 \in \X \logdot x_1 \in \Rb[x_2] \text{ and } y \sub \gamma x_1. \]
            It can be easily checked that $y \sub \gamma x_1 \miff x_1 \in \Inter y$, so we have
                \[ \neg \exists x_1 \in \X \logdot x_1 \in \Rb[x_2] \text{ and } x_1 \in \Inter y, \]
            implying that $\Rb[x_2] \inter \Inter y = \empt$.
            Since \axmref{tES.o.3} implies that $\Rb[x]$ is closed and $y$ consists of clopen, and, therefore, \textit{closed} subsets of $\X$, compactness implies that there is a finite subfamily $\bindset{K_i}{i=1}{n} \sub y$ such that
                \[ \Rb[x_2] \inter K_1 \inter \dots \inter K_n = \empt. \]
            Since $y \in \PrFilt(\clop{\X})$, we know that $y$ is closed under finite meets, implying $K_{\star} \in y$ where $K_{\star} = K_1 \inter \dots \inter K_n$.
            So we have $\Rb[x_2] \inter K_{\star} = \empt$, implying that $x_2 \notin \bDia K_{\star}$, implying that $\bDia K_{\star} \notin \gamma x_2$.
            But this gives us a contradiction, as it was assumed that $\gamma x_2 \Rb y$, but we have $K_{\star} \in y$ and $\bDia K_{\star} \notin \gamma x_2$.
    \end{prf}
    \begin{thm}\label{thm:duality}
        The functors $\spec{\ph} : \tHA \leftrightarrows \tES : \clop{\ph}$ are pseudo-inverse, implying that
            \[ \tHA \equiv \op{\tES}. \]
    \end{thm}
    \begin{prf}
        This follows directly from Lemmas \ref{lem:spec-well-defined-tha-obj} through \ref{lem:gamma-well-defined-nat-trans} as well as the above-mentioned fact that the naturality squares for $\pi$ and $\gamma$ commute.
    \end{prf}
    \par
    We now establish a congruence/filter/closed-upset correspondence for $\tHA$.
    This requires defining a subclass of filters called \mbox{« $\bDia$-filters »} and a subclass of subsets (on $\tTran$) called \mbox{« archival} \mbox{subsets »}.
    \begin{dfn}[$\bDia$-filter]
        Given $\A \in \tHA$ and $F \in \Filt(\A)$, we call $F$ a \textit{$\bDia$-filter} if for all $a,b \in \A$,
            \[ a \imp b \in F \quad \Mimp \quad \bDia a \imp \bDia b \in F. \]
        We denote the set of $\bDia$-filters on $\A$ by $\bDiaFilt(\A)$.
    \end{dfn}
    \begin{dfn}[Archival]\label{dfn:archival}
        Given $S \sub \X \in \tTran$, we say that $S$ is \textit{archival} if for all $x, z \in \X$,
            \[ x \notin S \contains z \text{ and } z \Rb x \quad \Mimp \quad {\Rb}[z] \inter \upset{x} \inter S \neq \empt. \]
        This is depicted as follows.
        
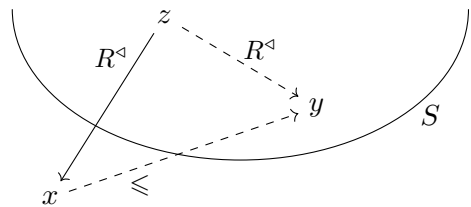
\begin{center}
        \begin{tikzpicture}
             \draw (0cm,2.5cm) arc(-180:0:3cm and 2cm);
             \node (S) at (5.5cm , 1.1cm) {$S$};
             \node (x) at (0.5cm , 0cm) {$x$};
             \node (y) at (4cm , 1.2cm) {$y$};
             \node (z) at (2cm , 2.4cm) {$z$};
             \draw (z) edge ["\small ${\Rb}$", ->, swap, near start, inner sep=1pt] (x);
             \draw (z) edge ["\small ${\Rb}$", ->, dashed, inner sep=1pt] (y);
             \draw (x) edge ["\small ${\leq}$", ->, dashed, swap, near start, inner sep=1pt] (y);
        \end{tikzpicture}
        \end{center}
        We denote the set of archival subsets of $\X$ by $\Arc(\X)$ and the set of archival \textit{upsets} of $\X$ by $\ArcUp(\X)$.
        If $\X$ has a topology defined on it, then we denote the the set of \textit{closed} archival upsets of $\X$ by $\ClArcUp(\X)$.
    \end{dfn}
    Note that since $x \neq z$, the fact that $z \Rb x$ is equivalent to the fact that $z > x$.
    This fact will be assumed at various points in the following proofs.
    \par
    We now prove the correspondence between congruences, $\bDia$-filters, and closed archival upsets.
    \begin{lem}\label{lem:congfilt-well-defined-tha}
        Given $\A \in \tHA$ and $\theta \in \Cong(\A)$, we have $\congfilt{\theta} \in \bDiaFilt(\A)$.
    \end{lem}
    \begin{prf}
        Since it is well known that $\congfilt{\theta} \in \Filt(\A)$, it remains only to show that $\congfilt{\theta}$ satisfies the additional $\bDia$-filter condition.
        Given $a \imp b \in \congfilt{\theta}$, we have $a \imp b \thetarel 1$, implying $a \meet b  = a \meet (a \imp b) \thetarel a \meet 1 = a$.
        Now since $a \meet b \thetarel a$, we have $\bDia (a \meet b) \thetarel \bDia a$, implying $\bDia a \thetarel \bDia (a \meet b)$.
        Also, since $a \meet b \thetarel a$, we have $b = (a \meet b) \join b \thetarel a \join b$, implying $\bDia b \thetarel \bDia (a \join b)$.
        Note that since $a \meet b \leq a \join b$, we have $\bDia(a \meet b) \leq \bDia(a \join b)$, implying $\bDia(a \meet b) \imp \bDia(a \join b) = 1$.
        Finally, since $\bDia a \thetarel \bDia (a \meet b)$ and $\bDia b \thetarel \bDia (a \join b)$, we have
            \[ \bDia a \imp \bDia b \thetarel \bDia (a \meet b) \imp \bDia (a \join b) = 1, \]
            implying $\bDia a \imp \bDia b \in \congfilt{\theta}$ as desired.
    \end{prf}
    \begin{lem}\label{lem:filtcong-well-defined-tha}
        Given $\A \in \tHA$ and $F \in \bDiaFilt(\A)$, we have $\filtcong{F} \in \Cong(\A)$.
    \end{lem}
    \begin{prf}
        Since it is well known that ${\filtcong{F}}$ is a congruence with respect to all non-modal operations, it remains only to show that ${\filtcong{F}}$ is a congruence over $\Box$ and $\bDia$.
        ~\\~\\
        $(\Box)$
            Given $a \filtcong{F} b$, we have $a \biimp b \in F$.
            Observe that we have the following inequalities.
            \[\begin{tikzcd}[ampersand replacement=\&,column sep=small,row sep=tiny]
                \& {a \imp b} \& {\Box(a \imp b)} \& {\Box a \imp \Box b} \\
                {a \biimp b} \\
                \& {b \imp a} \& {\Box (b \imp a)} \& {\Box b \imp \Box a}
                \arrow["{{\leq}}"{marking, allow upside down}, draw=none, from=1-2, to=1-3]
                \arrow["{{\leq}}"{marking, allow upside down}, draw=none, from=1-3, to=1-4]
                \arrow["{{\leq}}"{marking, allow upside down}, draw=none, from=2-1, to=1-2]
                \arrow["{{\leq}}"{marking, allow upside down}, draw=none, from=2-1, to=3-2]
                \arrow["{{\leq}}"{marking, allow upside down}, draw=none, from=3-2, to=3-3]
                \arrow["{{\leq}}"{marking, allow upside down}, draw=none, from=3-3, to=3-4]
            \end{tikzcd}\]
            This implies that $\Box a \imp \Box b, \Box b \imp \Box a \in F$, implying $\Box a \biimp \Box b \in F$, finally implying $\Box a \filtcong{F} \Box b$.
        ~\\~\\
        $(\bDia)$
            Given $a \filtcong{F} b$, we have $a \biimp b \in F$, implying $a \imp b, b \imp a \in F$, implying $\bDia a \imp \bDia b, \bDia b \imp \bDia a \in F$, implying $\bDia a \biimp \bDia b \in F$, finally implying $\bDia a \filtcong{F} \bDia b$.
    \end{prf}
    \begin{lem}\label{lem:filtclup-well-defined-tha}
        Given $\A \in \tHA$ and $F \in \bDiaFilt(\A)$, we have
            \[ \Inter \pi [F] \in \ClArcUp(\spec{\A}). \]
    \end{lem}
    \begin{prf}
        Since it is well known that $\Inter\pi[F] \in \ClUp(\spec{\A})$, it remains only to show that $\Inter\pi[F]$ is archival.
        Given $x, z \in \spec{\A}$ such that $x \sub z$ and $x \notin \Inter\pi[F] \contains z$, we show that there exists some $y \in \Rb[z] \inter \upset{x} \inter \Inter\pi[F]$.
        First note that since $x \notin \Inter\pi[F] \contains z$, we have $x \neq z$, implying that $x \Rf z$.
        Now consider the filter $\filter{F \union x}$ and the ideal $\bDia^{-1}[-z]$.
        Assume, toward a contradiction, that these sets are \textit{not} disjoint.
        Then there exists some $a \in F$ and some $b \in x$ such that $a \meet b \in \bDia^{-1}[-z]$, implying that $\bDia (a \meet b) \notin z$.
        Now since $b \in x$, we have $a \imp b \in x$, implying, since $x \Rf z$, that $\bDia(a \imp b) \in z$.
        Also, since $a \meet (a \imp b) = a \meet b$, we have, by weakening, that $a \meet (a \imp b) \leq a \meet b$, implying that $a \leq (a \imp b) \imp (a \meet b)$, implying that $(a \imp b) \imp (a \meet b) \in F$.
        Since $(a \imp b) \imp (a \meet b) \in F \in \bDiaFilt(\A)$, we have $\bDia (a \imp b) \imp \bDia(a \meet b) \in F \sub z$, implying, since $\bDia (a \imp b) \in z$, that $\bDia(a \meet b) \in z$,  contradicting our assumption that $\bDia(a \meet b) \notin z$.
        So the filter $\filter{F \union x}$ and the ideal $\bDia^{-1}[-z]$ are indeed disjoint, implying, by the Prime Filter Theorem, that there exists some $y \in \spec{\A}$ such that $\filter{F \union x} \sub y$ and $y \inter \bDia^{-1}[-z]= \empt$.
        Now we claim that $y \in \Inter\pi[F]$ and $x \sub y \Rf z$.
        To see that $y \in \Inter\pi[F]$ and $x \sub y$, simply observe that $F, x \sub \filter{F \union x} \sub y$.
        To see that the $z \Rb y$, take some $\bDia a \notin z$, implying that $a \in \bDia^{-1}[-z]$, finally implying $a \notin y$.
    \end{prf}
    \begin{lem}\label{lem:filtclup-inverse-well-defined-tha}
        Given $\A \in \tHA$ and $C \in \ClArcUp(\spec{\A})$, we have $\Inter C \in \bDiaFilt(\A)$.
    \end{lem}
    \begin{prf}
        Since it is well known that $\Inter C \in \Filt(\A)$, it remains only to show that $\Inter C$ satisfies the additional $\bDia$-filter condition.
        Assume, toward a contradiction, that $\Inter C$ does \textit{not} satisfy this condition, implying that there exist $a,b \in \A$ such that ${a \imp b \in \Inter C}$ but $\bDia a \imp \bDia b \notin \Inter C$.
        The latter of these facts implies that there exists some $w \in C$ such that $\bDia a \imp \bDia b \notin w$, implying that
            \[ w \notin \pi(\bDia a \imp \bDia b) = \bDia \pi a \imp \bDia \pi b = - \downset{(\bDia \pi a - \bDia \pi b)}, \]
        finally implying there exists some $z \in \spec{\A}$ such that $w \sub z$ and $z \in \bDia \pi a$ and $z \notin \bDia \pi b$.
        Note that since $w \sub z$ and $w \in C \in \Up(\spec{\A})$, we have $z \in C$.
        Also, since $z \in \bDia \pi a$, there exists some $x \in \pi a$ such that $z \Rb x$, but since $z \notin \bDia \pi b$ and $z \Rb x$, we have $x \notin \pi b$.
        Given $x \in \pi a$ and $x \notin \pi b$, we have $x \notin \pi a \imp \pi b$, implying $x \notin \pi (a \imp b)$, implying $a \imp b \notin x$.
        Since $a \imp b \in \Inter C$, we have $x \notin C$.
        So we have $z \Rb x$ and $x \notin C \contains z$, implying, since $C$ is archival, that there exists some $y \in C$ such that $z \Rb y$ and $x \sub y$.
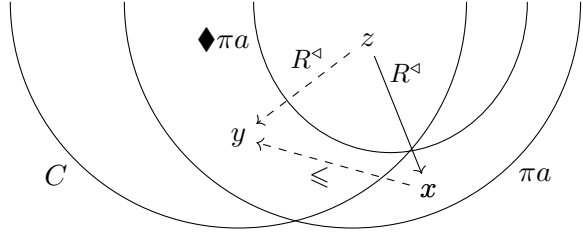
        \begin{center}
        \begin{tikzpicture}
             \draw (0,3cm) arc(-180:0:3cm and 3cm);
             \draw (1.5cm,3cm) arc(-180:0:3cm and 3cm);
             \draw (3.2cm,3cm) arc(-180:0:1.8cm and 2cm);
             \node (C) at (0.6cm,0.7cm) {$C$};
             \node (pa) at (6.9cm,0.7cm) {$\pi a$};
             \node (dpa) at (2.8cm,2.5cm) {$\bDia \pi a$};
             \node (C) at (5.5cm,0.5cm) {$x$};
             \node (x) at (5.5cm,0.5cm) {$x$};
             \node (y) at (3cm,1.2cm) {$y$};
             \node (z) at (4.7cm,2.5cm) {$z$};
             \draw (z) edge ["\small ${\Rb}$", ->, near start, inner sep=1pt] (x);
             \draw (z) edge ["\small ${\Rb}$", ->, swap, dashed, inner sep=1pt, near start] (y);
             \draw (x) edge ["\small ${\leq}$", ->, dashed, inner sep=1pt] (y);
        \end{tikzpicture}
        \end{center}
        Since $y \in C$, we have $a \imp b \in y$.
        Also, since $x \in \pi a$, we have $a \in x$, implying, since $x \sub y$, that $a \in y$, implying that $b \in y$.
        Finally, since $b \in y$ and $z \Rb y$, we have $\bDia b \in z$, contradicting the fact that $z \notin \bDia \pi b =  \pi \bDia b$.
    \end{prf}
    \begin{thm}\label{thm:cong-iso-bdiafilt-iso-clarcup}
        Given $\A \in \tHA$,
            \[ \tup{\Cong(\A), {\sub}} \isocat{\POS} \tup{\bDiaFilt(\A), {\sub}} \isocat{\POS} \tup{\ClArcUp(\spec{\A}), {\sup}}. \]
    \end{thm}
    \begin{prf}
        This follows directly from Lemmas \ref{lem:congfilt-well-defined-tha} through \ref{lem:filtclup-inverse-well-defined-tha} along with the fact that the relevant maps are monotone inverses of each other.
        \[\begin{tikzcd}[ampersand replacement=\&]
            {\tup{\Cong(\A), {\sub}}} \&\& {\tup{\bDiaFilt(\A), {\sub}}} \&\& {\tup{\ClArcUp(\spec{\A}), {\sup}}}
            \arrow["{\congfilt{\ph}}", shift left=2, tail, two heads, from=1-1, to=1-3]
            \arrow["{{\filtcong{\ph}}}", shift left=2, tail, two heads, from=1-3, to=1-1]
            \arrow["{\Inter\pi[\ph]}", shift left=2, tail, two heads, from=1-3, to=1-5]
            \arrow["\Inter", shift left=2, tail, two heads, from=1-5, to=1-3]
        \end{tikzcd}\]
    \end{prf}
    In the case that $\A \in \tHA$ is finite, we can characterise the $\bDia$-filters and, thereby, the congruences via \textit{elements} of $\A$.
    \begin{dfn}[$\bDia$-compatible]\label{dfn:bdia-compatible}
        Given $a \in \A \in \tHA$, we say that $a$ is \textit{$\bDia$-compatible} if for all $b \in \A$,
            \[ a \meet \bDia b \leq \bDia(a \meet b). \]
        We denote the set of $\bDia$-compatible elements of $\A$ by $\bDiaCom(\A)$.
        If $\bDiaCom(\A)$ has a second-greatest element, we call this element the \textit{$\bDia$-opremum}.
    \end{dfn}
    The nomenclature \mbox{« $\bDia$-compatible »} comes from \cite{CaiCig01}, where an operation $f: \A \to \A$ is said to be \textit{compatible} when, for all $a,b \in \A$,
        \[ a \meet fb = a \meet f(a \meet b) \]
    This can easily be shown to be equivalent to the condition given in Definition \ref{dfn:bdia-compatible} in our context.
    The operation $\bDia$ is \textit{not} compatible in general, but $\bDia$-compatible elements are exactly those elements that simulate compatibility and, for this reason, correspond to congruences in the finite case.
    \begin{prp}\label{prp:a-bdiacom-iff-upset-a-bdiafilt}
        Given $a \in \A \in \tHA$, we have
            \[ a \in \bDiaCom(\A) \quad \Miff \quad \upset{a} \in \bDiaFilt(\A). \]
    \end{prp}
    \begin{prf}
        $(\Rightarrow)$
            Clearly, $\upset{a} \in \Filt(\A)$, so we show that $\upset{a}$ satisfies the additional $\bDia$-filter condition.
            Given $b \imp c \in \upset{a}$, we have $a \leq b \imp c$, implying $a \meet b \leq c$, implying $\bDia(a \meet b) \leq \bDia c$.
            Since $a \in \bDiaCom(\A)$, we have $a \meet \bDia b \leq \bDia (a \meet b)$, so $a \meet \bDia b \leq \bDia c$, implying $a \leq \bDia b \imp \bDia c$, finally implying $\bDia b \imp \bDia c \in \upset{a}$.
        ~\\~\\
        $(\Leftarrow)$
            Arguing via the contrapositive, we assume $a \notin \bDiaCom(\A)$, implying there is some $b \in \A$ such that $a \meet \bDia b \notleq \bDia(a \meet b)$, implying that $a \notleq \bDia b \imp \bDia(a \meet b)$.
            Now since $a \meet b \leq a \meet b$, we have $a \leq b \imp (a \meet b)$.
            So we have $b \imp (a \meet b) \in \upset{a} $ but $\bDia b \imp \bDia(a \meet b)\notin \upset{a}$, implying that $\upset{a} \notin \bDiaFilt(\A)$.
    \end{prf}
    \begin{cor}\label{cor:bdiafilt-iso-bdiacom-fin}
        Given $\A \in \tHA_{\fin}$,
            \[ \tup{\bDiaFilt(\A), {\sub}} \isocat{\POS} \tup{\bDiaCom(\A), {\geq}}. \]
    \end{cor}
    \begin{prf}
        This follows directly from Proposition \ref{prp:a-bdiacom-iff-upset-a-bdiafilt} along with the fact that the relevant maps are monotone inverses of each other.
        \[\begin{tikzcd}[ampersand replacement=\&]
            {\tup{\bDiaFilt(\A), {\sub}}} \&\& {\tup{\bDiaCom(\A), {\geq}}}
            \arrow["\Meet", shift left=2, tail, two heads, from=1-1, to=1-3]
            \arrow["{\upset{}}", shift left=2, tail, two heads, from=1-3, to=1-1]
        \end{tikzcd}\]
    \end{prf}
\section{Notions of reachability}\label{sec:reachability}
    In  this section, we define and study two different notions of \mbox{« reachability » :} one on $\tES$ and the other on $\tTran_{\fin}$.
    These notions of reachability will be instrumental in dually characterising simple and subdirectly-irreducible temporal Heyting algebras.
    We show that in the finite case, these two notions coincide.
    \begin{dfn}
        Given $x,y \in \X \in \tES$, we say that $y$ \textit{is topo-reachable from} $x$ and write $x \tr y$ if $y$ is in every closed archival upset containing $x$, i.e.
            \[ y \in \Inter\compin{C}{\ClArcUp(\X)}{x \in C}. \]
        We call $x$ a \textit{topo-root} if it is a root with respect to the relation ${\tr}$, i.e.\ $\forall z \in \X \logdot x \tr z$.
        If such a point exists on $\X$, we say that $\X$ is \textit{topo-rooted}.
        We denote the set of topo-roots of $\X$ by $\ToRo(\X)$.
        Finally, we say that $\X$ is \textit{topo-connected} if it is connected with respect to the relation ${\tr}$, i.e.\ every point is a topo-root.
    \end{dfn}
    The relation ${\tr}$ is essentially analogous to the well-known \mbox{« specialisation} {order »} \cite[p.\ 37]{GehGoo24} except for the fact that $C$ ranges over the closed archival upsets of $\X$ instead of just the closed sets of $\X$.
    \par
    Here we prove a critical lemma for the characterisations in the next section.
    \begin{lem}\label{lem:toro-complement-arcup}
        Given $\X \in \tES$, we have $-\ToRo(\X) \in \ArcUp(\X)$.
    \end{lem}
    \begin{prf}
        We show that (1) $-\ToRo(\X)$ is an upset and (2) $-\ToRo(\X)$ is archival.
        ~\\~\\
        $(1)$
            We let $x \leq y$ and $x \in -\ToRo(\X)$ and show that $y \in -\ToRo(\X)$.
            Since $x \in -\ToRo(\X)$, there exists some $z \in \X$ such that $x \notbin{\tr} z$, implying that there exists some $C \in \ClArcUp(\X)$ such that $x \in C \notcontains z$.
            Since $x \in C \in \Up(\X)$ and $x \leq y$, we have $y \in C$.
            So we have $y \in C \in \ClArcUp(\X)$ and $z \notin C$, implying that $y \notbin{\tr} z$, further implying that $y$ is \textit{not} a topo-root, finally implying that $y \in -\ToRo(\X)$.
        ~\\~\\
        $(2)$
            We let $z \Rb x$ and $x \notin -\ToRo(\X) \contains z$ and find some $y \in \Rb[z] \inter \upset{x} \inter -\ToRo(\X)$.
            Since $z \notin \ToRo(\X)$, there exists some $w \in \X$ such that $z \notbin{\tr} w$, implying that there exists some $C \in \ClArcUp(\X)$ such that $z \in C \notcontains w$.
            Now it \textit{cannot} be the case that $x \in C$ because if it were, since $x \in \ToRo(\X)$ and, therefore, $x \tr w$, we would have $w \in C$, which we know not to be the case.
            So we have $z \Rb x$ and $x \notin C \contains z$, implying, since $C$ is archival, that there exists some $y \in \Rb[z] \inter \upset{x} \inter C$.
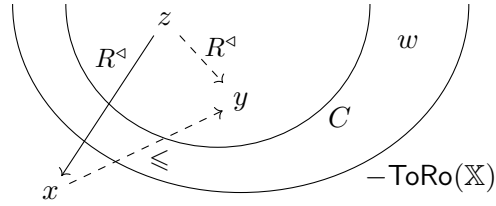
            \begin{center}
            \begin{tikzpicture}
                    \draw (0cm,2.5cm) arc(-180:0:3cm and 2.5cm);
                    \draw (.7cm,2.5cm) arc(-180:0:2cm and 1.9cm);
                    \node (C) at (4.3cm , 1cm) {$C$};
                    \node (T) at (5.5cm , 0.2cm) {$-\ToRo(\X)$};
                    \node (w) at (5.2cm , 2cm) {$w$};
                    \node (x) at (0.5cm , 0cm) {$x$};
                    \node (y) at (3cm , 1.2cm) {$y$};
                    \node (z) at (2cm , 2.3cm) {$z$};
                    \draw (z) edge ["\small ${\Rb}$", ->, swap, near start, inner sep=1pt] (x);
                    \draw (z) edge ["\small ${\Rb}$", ->, dashed, inner sep=1pt] (y);
                    \draw (x) edge ["\small ${\leq}$", ->, dashed, swap, inner sep=1pt] (y);
            \end{tikzpicture}
            \end{center}
            Since $y \in C \in \ClArcUp(\X)$ and $w \notin C$, we have $y \notbin{\tr} w$, implying that $y \notin \ToRo(\X)$.
            So we have found our $y \in \Rb[z] \inter \upset{x} \inter -\ToRo(\X)$.
    \end{prf}
    \par
    We now shift to the second notion of reachability, given solely in frame-theoretic terms and defined only on finite frames.
    \begin{dfn}
        Given $\X \in \tTran_{\fin}$, we define the following relation $B$.
            \[ x \Brel w \quad \Diff \quad w \leq x \text{ and } (w, x] \inter \Refl(\X) = \empt \]
        We also define the following relations for all $n \in \Nats$.
            \[ Z_0 := \Delta_{\X} \Quad Z_{n+1} := {Z_n}{;}{\Brel}{;}{\leq} \Quad Z := \Union_{m \in \Nats} Z_m \]
        Given $x,y \in \X$, we say that $y$ \textit{is $Z$-reachable from} $x$ if $x \Zrel y$.
        We call $x$ a \textit{$Z$-root} if it is a root with respect to the relation $Z$, i.e.\ $\forall z \in \X \logdot x \Zrel z$.
        If such a point exists on $\X$, we say that $\X$ is \textit{$Z$-rooted}.
        Finally, we say that $\X$ is \textit{$Z$-connected} if it is connected with respect to the relation $Z$.
    \end{dfn}
    The relation $B$ formalises the idea of moving \mbox{« backward »} until encountering either a reflexive or minimal point.
    The relation $Z$ then formalises the idea of \mbox{« zig-zagging »} down via $B$ and back up via ${\leq}$ a finite number of times.
    \begin{exm}
        Consider the following $\X \in \tTran_{\fin}$.
        \begin{multicols}{2}
        \noindent
        \[\begin{tikzcd}[ampersand replacement=\&,row sep=small]
            {} \& {z'} \& {} \& {} \\
            z \& {y'} \& {} \& {y'''} \\
            y \& {x'} \& {y''} \& {} \\
            x \& {w'} \& {} \& {} \\
            w \& {} \& {} \& {}
            \arrow[from=2-1, to=1-2]
            \arrow[from=3-1, to=2-1]
            \arrow[from=3-1, to=2-2]
            \arrow[from=3-3, to=2-2]
            \arrow[from=3-3, to=2-4]
            \arrow[from=4-1, to=3-1]
            \arrow[from=4-1, to=3-2]
            \arrow[from=4-1, to=4-1, loop, in=150, out=210, distance=5mm]
            \arrow[from=5-1, to=4-1]
            \arrow[from=5-1, to=4-2]
        \end{tikzcd}\]
        \columnbreak
        \noindent
        \begin{align*}
            Z_0[z] &= \set{z} \\
            B[Z_0[z]] &= \set{z, x, y} \\
            Z_1[z] &= \set{z, x, y, x', y', z'} \\
            B[Z_1[z]] &= \set{z, x, y, x', y', z', y''} \\
            Z_2[z] &= \set{z, x, y, x', y', z', y'', y'''} \\
            &\dots \\
            Z[z] &= \set{z, x, y, x', y', z', y'', y'''}
        \end{align*}
        \end{multicols}
        \noindent
        Here, everything is $Z$-reachable from $z$ except for $w$ and $w'$.
        This is because $B$ allows us to descend from $z$ to $x$, but does not allow us to descend to $w$ because $x \in (w,z] \inter \Refl(\X)$.
        Once we have descended as far as possible, we take upsets and repeat the process again.
    \end{exm}
    \par
    We prove several facts that will allow us to show that the two notions of reachability coincide in the finite case.
    \begin{prp}\label{prp:arc-characterisation-fin}
        Given $S \sub \X \in \tTran_{\fin}$, the subset $S$ is archival iff for all $x, z \in \X$,
            \[ z \Rb x \text{ and } x \notin S \contains z \quad \Mimp \quad \downset{z} \inter \upset{x} \inter \Refl(\X) \inter S \neq \empt. \]
        The consequent of the $\RHS$ of this biconditional can be depicted as follows.
        \begin{center}
        \begin{tikzpicture}
             \draw (0cm,2.5cm) arc(-180:0:3cm and 2cm);
             \node (S) at (5.5cm , 1.1cm) {$S$};
             \node (x) at (0.5cm , 0cm) {$x$};
             \node (y) at (4cm , 1.2cm) {$y$};
             \node (z) at (2cm , 2.4cm) {$z$};
             \draw (z) edge ["\small ${\Rb}$", ->, swap, near start, inner sep=1pt] (x);
             \draw (y) edge ["\small ${\leq}$", ->, dashed, swap, inner sep=1pt] (z);
             \draw (y) edge ["\small ${\Rf}$", ->, dashed, swap, out=330, in=60, looseness=6,  inner sep=1pt] (y);
             \draw (x) edge ["\small ${\leq}$", ->, dashed, swap, near start, inner sep=1pt] (y);
        \end{tikzpicture}
        \end{center}
    \end{prp}
    \begin{prf}
        $(\Rightarrow)$
            Given $z \Rb x$ and $x \notin S \contains z$, we know, since $S$ is archival, that $\Rb[z] \inter \upset{x} \inter S \neq \empt$.
            Since $\X$ is finite, every subset is well-founded, so consider some \textit{minimal} $y \in \Rb[z] \inter \upset{x}  \inter S$.
            Now since $y \in \Rb[z] \inter \upset{x}$, we have $x \leq y$ and $z \Rb y$, implying $y \Rf z$, further implying $x \leq y \leq z$, so we know that $y \in \downset{z} \inter \upset{x} \inter S$.
            So we need only show that $y \in \Refl(\X)$ for $y$ to witness the truth of the desired statement.
            Since we have $y \Rb x$ and $x \notin S \contains y$ and $S$ is archival, there exists some $w \in \Rb[y] \inter \upset{x} \inter S$.
            Now since $w \in \Rb[y]$, we have $w \Rf y$, implying $w \leq y$.
            But since $y$ was assumed to be minimal, it must be the case that $y = w \Rf y$, implying that $y \in \Refl(\X)$ as desired.
        ~\\~\\
        $(\Leftarrow)$
            We let $\downset{z} \inter \upset{x} \inter \Refl(\X) \inter S \neq \empt$ and show that $\Rb[z] \inter \upset{x} \inter S \neq \empt$.
            Take $y \in \Refl(\X) \inter S$ such that $x \leq y \leq z$.
            Since $y \leq y \Rf y \leq z$, the \mbox{« mix »} condition implies that $y \Rf z$, implying that $y \in \Rb[z] \inter \upset{x} \inter S$, and, therefore, that $S$ is archival.
    \end{prf}
    \begin{lem}\label{lem:arc-closed-under-b}
        Archival subsets are closed under $B$ (on \textit{finite} temporal transits, the only frames on which $B$ is defined).
    \end{lem}
    \begin{prf}
        Let $x \in S \in \Arc(\X)$ and $x \Brel w$.
        We show that $w \in S$ via induction on the length of the longest chain in $(w, x]$, denoted by $n$.
        This chain must be finite because $B$ is only defined on finite frames.
        ~\\~\\
        $(n = 0)$
            This means that $w = x$, implying, since $x \in S$, that $w \in S$.
        ~\\~\\
        $(n = m + 1)$
            This means that we have a chain of unique elements $\bindset{y_i}{i=1}{m+1} \sub (w, x]$ such that
                \[ w < y_1 < \dots < y_m < y_{m+1} = x. \]
            Now since $x \Brel w$, we have $(w, x] \inter \Refl(\X) = \empt$, implying $(y_1, x] \inter \Refl(\X) = \empt$.
            Combining this with the fact that $y_1 \leq x$, and that the length of the longest chain in $(y_1, x]$ is $\leq m$, we can apply our induction hypothesis, implying that $y_1 \in S$.
            Now assume, toward a contradiction, that $w$ is \textit{not} in $S$.
            This implies that $w \leq y_1$ and $w \notin S \contains y_1$, implying, since $S$ is archival, that there exists some $y \in S$ such that $w \leq y$ and $y_1 \Rb y$.
            Now we cannot have $y = y_1$ because then $y_1 \Rb y_1$, contradicting $y_1 \notin \Refl(\X)$.
            We also cannot have $w = y$ because $w \notin S \contains y$.
            So it must be the case that $w < y < y_1$.
            But observe that this contradicts the fact that $\bindset{y_i}{i=1}{m+1}$ was taken to be the longest chain in $(w, x]$.
    \end{prf}
    \begin{lem}\label{lem:arcup-closed-under-z}
        Archival upsets are closed under $Z$.
    \end{lem}
    \begin{prf}
        Let $w \in U \in \ArcUp(\X)$ and $w \Zrel z$.
        Since $w \Zrel z$, we have $w \mathbin{Z_n} z$ for some $n \in \Nats$.
        We proceed via induction on $n$.
        ~\\~\\
        $(n=0)$
            Given $w \Zrel_{0} z$, we have $w = z$, implying that $z \in U$.
        ~\\~\\
        $(n=m+1)$
            Given $w \mathbin{Z_{m+1}} z$, we have $w \mathbin{Z_m} x \Brel y \leq z$ for some $x,y \in \X$.
            By induction hypothesis, we have $x \in U$.
            Since $x \Brel y$ and $x \in U \in \Arc(\X)$, we know, by Lemma \ref{lem:arc-closed-under-b}, that $y \in U$.
            And since $y \leq z$ and $y \in U \in \Up(\X)$, we can conclude that $z \in U$.
    \end{prf}
    \begin{lem}\label{lem:z-image-arcup-fin}
        Given $S \sub \X \in \tTran_{\fin}$, we have $Z[S] \in \ArcUp(\X)$.
    \end{lem}
    \begin{prf}
        We show that (1) $Z[S]$ is an upset and (2) $Z[S]$ is archival.
        ~\\~\\
        $(1)$
            We let $y \in Z[S] $ and $y \leq z$ and show that $z \in Z[S]$.
            Now since $y \in Z[S]$, we have $x \in S$ such that $x \mathbin{Z_n} y$ for some $n \in \Nats$.
            This implies that we have $x \mathbin{Z_n} y \Brel y \leq z$, implying that we have $x \mathbin{Z_{n+1}} z$ as desired.
        ~\\~\\
        $(2)$
            Now, to see that $Z[S]$ is archival, we let $x \leq z$ and $x \notin Z[S] \contains z$.
            Recall that Proposition \ref{prp:arc-characterisation-fin} implies that it suffices to show that
                \[ \downset{z} \inter \upset{x} \inter \Refl(\X) \inter Z[S] \neq \empt. \]
            Now since $x \notin Z[S]$, it cannot be the case that $z \Brel x$ as we would then have $z \Brel x \leq x$, implying $x \in Z[S]$.
            Since $x \leq z$ and $z \notbin{B} x$, it must be the case that $(x, z] \inter \Refl(\X) \neq \empt$.
            Since $\X$ is finite, let $y$ be some \textit{maximal} point in $(x, z] \inter \Refl(\X)$.
            This implies that $y \leq z$ and $(y, z] \inter \Refl(\X) = \empt$, implying that $z \Brel y$, implying, by Lemma \ref{lem:arc-closed-under-b}, that $y \in Z[S]$.
            So we have found our point
                \[ y \in \downset{z} \inter \upset{x} \inter \Refl(\X) \inter Z[S]. \]
    \end{prf}
    \par
    We now show that the two above-defined notions of reachability coincide in the finite case.
    \begin{prp}\label{cor:coincidence-reachability}
        Given $x, y \in \X \in \tES_{\fin}$, $y$ is topo-reachable from $x$ iff $y$ is $Z$-reachable from $x$ (on the underlying frame of $\X$).
    \end{prp}
    \begin{prf}
        $(\Rightarrow)$
            Lemma \ref{lem:z-image-arcup-fin} implies $Z[x] \in \ArcUp(\X)$, implying, since the topology is discrete, that $Z[x] \in \ClArcUp(\X)$.
            Since $x \in Z[x] \in \ClArcUp(\X)$ and $x \tr y$, we have $y \in Z[x]$, implying that $y$ is $Z$-reachable from $x$.
        ~\\~\\
        $(\Leftarrow)$
            We let $x \in C \in \ClArcUp(\X)$ and show that $y \in C$.
            Since $x \Zrel y$ and $x \in C \in \ArcUp(\X)$, Lemma \ref{lem:arcup-closed-under-z} implies that $y \in C$ as desired.
    \end{prf}
    \begin{cor}\label{cor:coincidence-rooted-connected}
        Given $x\in \X \in \tES_{\fin}$, $x$ is a topo-root iff $x$ is a $Z$-root (on the underlying frame of $\X$).
        Furthermore, $\X$ is topo-connected iff its underlying frame is $Z$-connected.
    \end{cor}
\section{Characterisations}\label{sec:characterisations}
    In this section, we characterise simple and subdirectly-irreducible temporal Heyting algebras lattice-theoretically and order-topologically.
    In the finite cases, we do the same element-wise and frame-theoretically.
    Note that we make liberal use of the well-known characterisation of subdirectly-irreducible algebras as exactly those algebras having a « monolith » (a second-least congruence, equiv.\ a least non-$\Delta$ congruence) \cite[Theorem 3.23]{Ber12}.
    \par
    We first characterise \textit{simple} algebras lattice-theoretically and order-topologically.
    \begin{thm}\label{thm:simple-characterisation-general}
        Given $\A \in \tHA$, the following are equivalent.
        \begin{enumerate}
            \item $\A$ is simple
            \item $\bDiaFilt(\A) = \set{\set{1}, \A}$
            \item $\spec{\A}$ is topo-connected
        \end{enumerate}
    \end{thm}
    \begin{prf}
        Let $\X := \spec{\A}$ and recall Theorem \ref{thm:cong-iso-bdiafilt-iso-clarcup} :
            \[ \tup{\Cong(\A), {\sub}} \isocat{\POS} \tup{\bDiaFilt(\A), {\sub}} \isocat{\POS} \tup{\ClArcUp(\spec{\A}), {\sup}}. \]
        $(1 \Leftrightarrow 2)$
            This follows directly from Theorem \ref{thm:cong-iso-bdiafilt-iso-clarcup}.
        ~\\~\\
        $(2 \Rightarrow 3)$
            Here we take some $x, y \in \X$ and show that $y$ is topo-reachable from $x$.
            Since $\bDiaFilt(\A) = \set{\set{1}, \A}$, we have
                \[ \ClArcUp(\A) = \Set{\Inter\pi[\set{1}], \Inter\pi[\A]} = \Set{\X, \empt}. \]
            This implies that the only closed archival upset containing $x$ is $\X$, so
                \[ y \in \X = \Inter \set{\X} = \Inter \compin{C}{\ClArcUp(\X)}{x \in C}, \]
            implying that $y$ is topo-reachable from $x$ as desired.
        $(2 \Leftarrow 3)$
            Arguing via the contrapositive, suppose we have $\bDiaFilt(\A) \neq \set{\set{1}, \A}$, implying that there is some $F \in \bDiaFilt(\A)$ such that $\set{1} \neq F \neq \A$.
            This implies that $\Inter\pi[F] \in \ClArcUp(\X)$ and
                \[ \X = \Inter\pi[\set{1}] \neq \Inter\pi[F] \neq \Inter\pi[\A] = \empt. \]
            Since $\X \neq \Inter\pi[F] \neq \empt$, we have some $x,y \in \X$ such that $x \in \Inter\pi[F] \notcontains y$.
            Since $x \in \Inter\pi[F] \in \ClArcUp(\X)$ and $y \notin \Inter\pi[F]$, $y$ is \textit{not} topo-reachable from $x$, implying that $\X$ is \textit{not} topo-connected.
    \end{prf}
    We now characterise \textit{finite} simple algebras element-wise and frame-theoretically.
    \begin{cor}\label{thm:cor-characterisation-fin}
        Given $\A \in \tHA_{\fin}$, the following are equivalent.
        \begin{enumerate}
            \item $\A$ is simple
            \item $\bDiaCom(\A) = \set{1, 0}$
            \item $\fram{\A}$ is $Z$-connected
        \end{enumerate}
    \end{cor}
    \begin{prf}
        The proof is summarised by the following diagram.
        \[\begin{tikzcd}
            & {(1)} \\
            \\
            {\bDiaFilt(\A) = \set{\set{1}, \A}} && {\spec{\A} \text{ topo-connected}} \\
            \\
            {(2)} && {(3)}
            \arrow["{\text{Thm.\ \ref{thm:simple-characterisation-general}}}"{description}, between={0.1}{0.9}, Rightarrow, 2tail reversed, from=3-1, to=1-2]
            \arrow["{\text{Thm.\ \ref{thm:simple-characterisation-general}}}"{description}, between={0.1}{0.9}, Rightarrow, 2tail reversed, from=3-1, to=3-3]
            \arrow["{\text{Cor.\ \ref{cor:bdiafilt-iso-bdiacom-fin}}}"{description}, between={0.1}{0.9}, Rightarrow, 2tail reversed, from=3-1, to=5-1]
            \arrow["{\text{Thm.\ \ref{thm:simple-characterisation-general}}}"{description}, between={0.1}{0.9}, Rightarrow, 2tail reversed, from=3-3, to=1-2]
            \arrow["{\text{Cor.\ \ref{cor:coincidence-rooted-connected}}}"{description}, between={0.1}{0.9}, Rightarrow, 2tail reversed, from=3-3, to=5-3]
        \end{tikzcd}\]
    \end{prf}
    \par
    Next we characterise \textit{subdirectly-irreducible} algebras lattice-theoretically and order-topologically.
    \begin{thm}\label{thm:si-characterisation-general}
        Given $\A \in \tHA$, the following are equivalent.
        \begin{enumerate}
            \item $\A$ is subdirectly-irreducible
            \item $\bDiaFilt(\A)$ has a second-least element
            \item $\ToRo(\spec{\A})$ is non-empty and open
        \end{enumerate}
    \end{thm}
    \begin{prf}
        Let $\X := \spec{\A}$ and recall Theorem \ref{thm:cong-iso-bdiafilt-iso-clarcup} :
            \[ \tup{\Cong(\A), {\sub}} \isocat{\POS} \tup{\bDiaFilt(\A), {\sub}} \isocat{\POS} \tup{\ClArcUp(\spec{\A}), {\sup}}. \]
        $(1 \Leftrightarrow 2)$
            This follows directly from Theorem \ref{thm:cong-iso-bdiafilt-iso-clarcup}.
        ~\\~\\
        $(2 \Rightarrow 3)$
            Let $F$ be the second-least $\bDia$-filter, implying that $\set{1} \neq F$ and that $\Inter\pi[F] \neq \X$ is the second-greatest closed archival upset.
            We claim that
                \[ \Inter\pi[F] = -\ToRo(\X). \]
            Note that if this is the case, then we'll have $-\ToRo(\X) \in \ClArcUp(\X)$, implying that $\ToRo(\X)$ is open.
            Also, since $\Inter\pi[F] \neq \X$, we'll have $-\ToRo(\X) \neq \X$ and, therefore, $\empt \neq \ToRo(\X)$, implying that $\ToRo(\X)$ is non-empty as desired.
            So we show that $\Inter\pi[F] = -\ToRo(\X)$.
            $(\sub)$
                Since $\Inter\pi[F] \neq \X$, we have some $y \notin \Inter\pi[F]$.
                Given $x \in \Inter\pi[F]$, we have $x \in \Inter\pi[F] \in \ClArcUp(\X)$ and $y \notin \Inter\pi[F]$, implying that $y$ is \textit{not} topo-reachable from $x$, implying that $x \in -\ToRo(\X)$.
            $(\sup)$
                Arguing via the contrapositive, let $x \notin \Inter\pi[F]$.
                Since $\Inter\pi[F]$ is the second-largest closed archival upset, $x \notin \Inter\pi[F]$ implies that the only closed archival upset containing $x$ is $\X$.
                Given an arbitrary $y \in \X$, we then have
                    \[ y \in \X = \Inter \set{\X} = \Inter \compin{C}{\ClArcUp(\X)}{x \in C}, \]
                implying that $x$ is a topo-root and, therefore, $x \notin -\ToRo(\X)$.
        ~\\~\\
        $(2 \Leftarrow 3)$
            Assume toward a contradiction, that $\ToRo(\X)$ is non-empty and open, but $\bDiaFilt(\A)$ does \textit{not} have a second-least element, implying that
                \[ \forall F \neq \set{1} \logdot \exists F' \neq \set{1} \logdot F \notsub F' \]
            (where $F, F'$ are taken to range over $\bDiaFilt(\A)$).
            This implies that
                \[ \forall C \neq \X \logdot \exists C' \neq \X \logdot C' \notsub C \]
            (where $C, C'$ are taken to range over $\ClArcUp(\X)$).
            But observe that if\\$\ToRo(\X)$ is non-empty and open, then $-\ToRo(\X)$ is non-total and closed, implying, by Lemma \ref{lem:toro-complement-arcup}, that $-\ToRo(\X) \in \ClArcUp(\X)$.
            So there exists some $C' \in \ClArcUp(\X)$ such that $\X \neq C' \notsub -\ToRo(\X)$, implying that we have some $x \in \X$ such that $C' \contains x \notin -\ToRo(\X)$, and, therefore, that $x$ is a topo-root.
            But since $\X \neq C'$, there exists some $y \notin C'$, implying, since $x \in C' \in \ClArcUp(\X)$, that $x \notbin{\tr} y$, contradicting the fact that $x$ is a topo-root.
    \end{prf}
    Finally, we characterise \textit{finite} subdirectly-irreducible algebras element-wise and frame-theoretically.
    \begin{cor}\label{cor:si-characterisation-fin}
        Given $\A \in \tHA_{\fin}$, the following are equivalent.
        \begin{enumerate}
            \item $\A$ is subdirectly-irreducible
            \item $\A$ has a $\bDia$-opremum
            \item $\fram{\A}$ is $Z$-rooted
        \end{enumerate}
    \end{cor}
    \begin{prf}
        The proof is summarised by the following diagram.
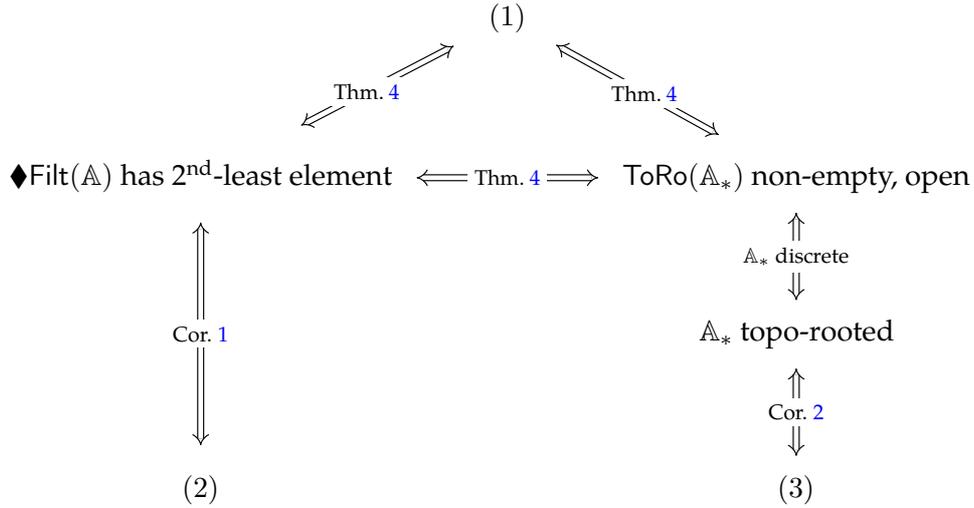
        \[\begin{tikzcd}
            & {(1)} \\
            \\
            \begin{array}{c} \bDiaFilt(\A) \text{ has a}\\\text{second-least element} \end{array} && \begin{array}{c} \ToRo(\spec{\A}) \text{ is}\\\text{non-empty and open} \end{array} \\
            \\
            && {\spec{\A} \text{ topo-rooted}} \\
            \\
            {(2)} && {(3)}
            \arrow["{\text{Thm.\ \ref{thm:si-characterisation-general}}}"{description}, between={0.1}{0.9}, Rightarrow, 2tail reversed, from=3-1, to=1-2]
            \arrow["{\text{Thm.\ \ref{thm:si-characterisation-general}}}"{description}, between={0.1}{0.9}, Rightarrow, 2tail reversed, from=3-1, to=3-3]
            \arrow["{\text{Cor.\ \ref{cor:bdiafilt-iso-bdiacom-fin}}}"{description}, between={0.1}{0.9}, Rightarrow, 2tail reversed, from=3-1, to=7-1]
            \arrow["{\text{Thm.\ \ref{thm:si-characterisation-general}}}"{description}, between={0.1}{0.9}, Rightarrow, 2tail reversed, from=3-3, to=1-2]
            \arrow["{\spec{\A} \text{ discrete}}"{description}, between={0.1}{0.9}, Rightarrow, 2tail reversed, from=5-3, to=3-3]
            \arrow["{\text{Cor.\ \ref{cor:coincidence-rooted-connected}}}"{description}, between={0.1}{0.9}, Rightarrow, 2tail reversed, from=7-3, to=5-3]
        \end{tikzcd}\]
    \end{prf}
    \par
    To conclude this section, we give examples of finite subdirectly-irreducible and simple temporal Heyting algebras respectively, along with their dual frames, to visually illustrate the characterisations given above.
    On the algebras, the ${\Box}$ operation is represented by arrows with an empty dot and the ${\bDia}$ operation by arrows with a filled dot.
    Fixpoints are omitted.
    On their dual frames, the same minimal method of representing temporal transits is employed as in Example \ref{exm:transit}.
    \begin{exm}\label{exm:fin-si}
        The following is an example of a finite \textit{subdirectly-irreducible} temporal Heyting algebra $\A$ along with its $Z$-rooted dual frame $\fram{\A}$.
        \[\begin{tikzcd}[sep=scriptsize]
            && 1 &&&&& \\
            & e && f && {\upset{a}} && {\upset{b}} \\
            c && d \\
            & a && b & {\upset{c}} && {\upset{f}} \\
            && 0
            \arrow[no head, from=2-2, to=1-3]
            \arrow[no head, from=2-4, to=1-3]
            \arrow[from=2-8, to=2-8, loop, in=75, out=15, distance=5mm]
            \arrow[no head, from=3-1, to=2-2]
            \arrow[no head, from=3-3, to=2-2]
            \arrow[no head, from=3-3, to=2-4]
            \arrow["\bullet"{marking}, curve={height=12pt}, from=3-3, to=4-4]
            \arrow[no head, from=4-2, to=3-1]
            \arrow[no head, from=4-2, to=3-3]
            \arrow["\bullet"{marking}, curve={height=12pt}, from=4-2, to=5-3]
            \arrow[no head, from=4-4, to=3-3]
            \arrow["\bullet"{marking, text=\pgfkeysvalueof{/tikz/commutative diagrams/background color}}, "\circ"{marking}, curve={height=12pt}, from=4-4, to=3-3]
            \arrow[from=4-5, to=2-6]
            \arrow[from=4-5, to=4-5, loop, in=345, out=285, distance=5mm]
            \arrow[from=4-7, to=2-6]
            \arrow[from=4-7, to=2-8]
            \arrow[from=4-7, to=4-7, loop, in=345, out=285, distance=5mm]
            \arrow[no head, from=5-3, to=4-2]
            \arrow["\bullet"{marking, text=\pgfkeysvalueof{/tikz/commutative diagrams/background color}}, "\circ"{marking}, curve={height=12pt}, from=5-3, to=4-2]
            \arrow[no head, from=5-3, to=4-4]
        \end{tikzcd}\]
        It is of interest to note that the $\HA$-reduct of $\A$ is \textit{not} subdirectly-irreducible as it does not have a second-least filter and the dual frame is not rooted with respect to ${\leq}$.
        To analyse $\A$ from the point of view of Corollary \ref{cor:si-characterisation-fin}, the $\bDia$-opremum is $b$ and the points $\upset{c}$, $\upset{a}$, and $\upset{f}$ are all $Z$-roots on $\fram{\A}$.
        Observe that $\A$ is \textit{not} simple as it has a non-trivial $\bDia$-compatible element ($b$) and $\spec{\A}$ has a non-$Z$-root ($\upset{b}$).
    \end{exm}
    \begin{exm}\label{exm:fin-simple}
        The following is an example of a finite \textit{simple} temporal Heyting algebra $\B$ along with its $Z$-connected dual frame $\fram{\B}$.
        \[\begin{tikzcd}[sep=scriptsize]
            && 1 &&&&& \\
            & e && f && {\upset{a}} && {\upset{b}} \\
            c && d \\
            & a && b & {\upset{c}} && {\upset{f}} \\
            && 0
            \arrow[no head, from=2-2, to=1-3]
            \arrow["\bullet"{marking}, curve={height=12pt}, from=2-2, to=3-1]
            \arrow[no head, from=2-4, to=1-3]
            \arrow[no head, from=3-1, to=2-2]
            \arrow["\bullet"{marking, text=\pgfkeysvalueof{/tikz/commutative diagrams/background color}}, "\circ"{marking}, curve={height=12pt}, from=3-1, to=2-2]
            \arrow[no head, from=3-3, to=2-2]
            \arrow[no head, from=3-3, to=2-4]
            \arrow["\bullet"{marking}, curve={height=12pt}, from=3-3, to=5-3]
            \arrow[no head, from=4-2, to=3-1]
            \arrow[no head, from=4-2, to=3-3]
            \arrow["\bullet"{marking, text=\pgfkeysvalueof{/tikz/commutative diagrams/background color}}, "\circ"{marking}, curve={height=-12pt}, from=4-2, to=3-3]
            \arrow["\bullet"{marking}, curve={height=12pt}, from=4-2, to=5-3]
            \arrow[no head, from=4-4, to=3-3]
            \arrow["\bullet"{marking, text=\pgfkeysvalueof{/tikz/commutative diagrams/background color}}, "\circ"{marking}, curve={height=12pt}, from=4-4, to=3-3]
            \arrow["\bullet"{marking}, curve={height=-12pt}, from=4-4, to=5-3]
            \arrow[from=4-5, to=2-6]
            \arrow[from=4-5, to=4-5, loop, in=345, out=285, distance=5mm]
            \arrow[from=4-7, to=2-6]
            \arrow[from=4-7, to=2-8]
            \arrow[from=4-7, to=4-7, loop, in=345, out=285, distance=5mm]
            \arrow["\bullet"{marking, text=\pgfkeysvalueof{/tikz/commutative diagrams/background color}}, "\circ"{marking}, curve={height=12pt}, from=5-3, to=3-3]
            \arrow[no head, from=5-3, to=4-2]
            \arrow[no head, from=5-3, to=4-4]
        \end{tikzcd}\]
        Observe that $\B$ is now simple (compared to $\A$ in Example \ref{exm:fin-si}) as $b$ is no longer $\bDia$-compatible :
            \[ b \meet \bDia f = b \meet f = b \nleq 0 = \bDia b = \bDia(b \meet f). \]
        This contrasts with Example \ref{exm:fin-si}, where we had
            \[ b \meet \bDia f = b \meet f = b \leq b = \bDia b = \bDia(b \meet f). \]
        On the dual frame $\fram{\B}$, we now have $Z$-connectedness as the point $\upset{b}$ is no longer ${\Rf}$-reflexive, allowing it to $Z$-reach all other points (whereas in Example \ref{exm:fin-si} it could only $Z$-reach itself).
    \end{exm}
\section{Applications to tHC}\label{sec:applications}
    In this section, we apply several results proven above to study the temporal Heyting calculus.
    \par
    We use Theorem \ref{thm:duality} to prove the relational and algebraic FMP for the logic.
    Duality allows us to work with the more friendly \textit{relational} filtration and transfer finiteness back to algebras.
    Indeed, while algebraic filtrations of the intuitionistic ${\imp}$ are well-understood \cite{BezBez16}, algebraic filtrations of \textit{modal operators} on Heyting algebras do not seem to have been covered anywhere in the literature and appear, \textit{prima facie}, to be quite unwieldy.
    \par
    We then use the algebraic FMP in conjunction with Corollary \ref{cor:si-characterisation-fin} and some well-known results from universal algebra to prove a relational completeness result for the temporal Heyting calculus that combines finiteness and the frame property dual to subdirect-irreducibility.
    This is analogous to the fact that $\IPC$ is complete with respect to the class of finite rooted posets \cite[Theorem 6.12]{TroDal88} as rootedness corresponds to subdirect-irreducibility for $\HA$ in the finite case \cite[Theorem 2.3.16]{Bez06}.
    \par
    Roughly, we follow the following path to $\mathbf{Z}$, the class of finite $Z$-rooted temporal transits, using duality to avoid the less-friendly routes involving algebraic filtration and \mbox{« unravelling »} \cite[\textsection4.5]{BlaRijVen01}.
    \[\begin{tikzcd}
        & \tHC \\
        \\
        \\
        \tHA && \tES \\
        \\
        {\tHA_{\textsf{fin}}} && {\tTran_{\textsf{fin}}} \\
        \\
        {\tHA_{\textsf{fsi}}} && {\mathbf{Z}}
        \arrow["{\text{completeness}}", no head, from=4-1, to=1-2]
        \arrow["{\text{duality}}", no head, from=4-1, to=4-3]
        \arrow["{\text{algebraic filtration}}"', dashed, from=4-1, to=6-1]
        \arrow["{\text{completeness}}"', no head, from=4-3, to=1-2]
        \arrow["{\text{relational filtration}}", Rightarrow, from=4-3, to=6-3]
        \arrow["{\text{Birkoff's Theorem}}"', Rightarrow, from=6-1, to=8-1]
        \arrow["{\text{duality}}"', Rightarrow, from=6-3, to=6-1]
        \arrow["{\text{unravelling}}", dashed, from=6-3, to=8-3]
        \arrow["{\text{+ characterisations}}"', draw=none, from=8-1, to=8-3]
        \arrow["{\text{duality}}", Rightarrow, from=8-1, to=8-3]
    \end{tikzcd}\]
	\vspace{3pt}
    \par
    Letting $\Prop$ be a fixed set of propositional variables, we define the following languages.
        \[ \Lm := p \in \Prop \mid \phi \land \phi \mid \phi \lor \phi \mid \phi \imp \phi \mid \Box \phi \mid \bot \mid \top \quad\quad \Lt := \phi \in \Lm \mid \bDia \phi \]
    Here we recall the axiomatisations of the \textit{modalized Heyting calculus} and \textit{temporal Heyting calculus} as given in \cite{Esa06}.
    \begin{dfn}
        The \textit{modalized Heyting calculus}, denoted by $\mHC$, is the smallest subset of $\Lm$ that contains $\IPC$, the following axioms, and is closed under modus ponens and uniform substitution.
            \[ \Box (p \imp q) \imp (\Box p \imp \Box q) \Quad p \imp \Box p \Quad \Box p \imp (q \lor (q \imp p)) \]
    \end{dfn}
    We define the rule \axm{PD} as follows for formulas $\phi, \chi \in \Lt$ : \raisebox{-9pt}{$\prftree{\phi \imp \chi}{\bDia \phi \imp \bDia \chi}$} ~.
    \label{PD}
    \begin{dfn}
        The \textit{temporal Heyting calculus}, denoted by $\tHC$, is the smallest subset of $\Lt$ that contains $\mHC$, the following axioms, and is closed under modus ponens, uniform substitution, and \axmref{PD}.
            \[ \bDia (p \lor q) \imp (\bDia p \lor \bDia q) \Quad \bDia \bot \imp \bot \Quad p \imp \Box\bDia p \Quad \bDia \Box p \imp p \]
    \end{dfn}
    \par
    It can be shown, given the standard \textit{algebraic} semantics, that we have $\tHC \soco \tHA$ (read « $\tHC$ is sound and complete with respect to $\tHA$ »).
    Completeness is achieved via the well-known \mbox{« Lindenbaum-Tarski} \mbox{process »} \cite{FonJanPig03}.
    Given this fact, an \textit{algebraic model} (of $\tHC$) is a tuple $\tup{\A, \nu}$ such that $\A \in \tHA$ and $\nu$ is an algebraic valuation, i.e.\ a $\tHA$-morphism $\Term \to \A$, where $\Term$ is the algebra of terms on $\Prop$ constructed with the language $\Lt$.
    \par
    It can also be shown, given the standard \textit{relational} semantics (reading $\Box$ as a forward-looking universal and $\bDia$ as a backward-looking existential), that we have $\tHC \soco \tTran$.
    Completeness can be achieved either via the method of canonical models \cite[\textsection4.2]{BlaRijVen01} or via Esakia duality.
    Given this fact, a \textit{relational model} (of $\tHC$) is a tuple $\tup{\X, \nu}$ such that $\X \in \tTran$ and $\nu$ is an (intuitionistic) relational valuation, i.e.\ a map $\Prop \to \Up(\X)$ .
    \par
    For details on these algebraic and relational completeness results, see \cite[\textsection3.1]{Alv24} and \cite[\textsection6.2]{Alv24} respectively.
    \par
    We now define a means of turning an algebraic model into a relational model.
    We define a slight abuse-of-notation, extending the functor $\spec{\ph} : \tHA \to \tES$ to map between algebraic and relational \textit{models} of $\tHC$.
    \begin{dfn}\label{dfn:dual-relational-model}
        Given an algebraic model $\M := \tup{\A, \nu}$, the \textit{dual relational model} of $\M$ is defined as
            \[ \spec{\M} := \tup{\fram{\A}, \pi \after \nu}. \]
    \end{dfn}
    Note that since $\pi : \A \to \ClopUp(\spec{\A})$, we have $\pi \after \nu : \Term \to \ClopUp(\spec{\A})$, implying that $\pi \after \nu$, when restricted to $\Prop \sub \Term$, is a relational valuation with respect to the poset relation ${\leq}$, as it maps onto $\ClopUp(\spec{\A}) \sub \Up(\fram{\A})$.
    Note, also, that it can easily be shown that any valuation with domain $\Prop$ uniquely determines a valuation with domain $\Term$ and vice-versa.
    \par
    We now relate the theories of an algebraic model and its dual relational model.
    \begin{lem}[Algebraic truth lemma]\label{lem:truth-lemma}
        Given an algebraic model $\M$, a formula $\phi \in \Lt$, and a prime filter $x \in \spec{\M}$,
            \[ \nu\phi \in x \quad \Miff \quad \tup{\spec{\M}, x} \models \phi. \]
    \end{lem}
    \begin{prf}
        We proceed via induction on the shape of $\phi$, but show only the case where $\phi$ is of the form $\bDia\chi$ as the other cases are trivial, well-known, or follow by symmetrical arguments.
        Note that this proof makes use of the fact that $\pi$ and $\nu$ are both $\tHA$-morphisms, and, therefore, homomorphic over ${\bDia}$.
        \begin{align*}
            \nu \bDia \chi \in x \quad &\Miff \quad x \in \pi \nu \bDia \chi \\
            &\Miff \quad x \in \bDia \pi \nu  \chi \\
            &\Miff \quad \Rb[x] \inter \pi \nu \chi \neq \empt \\
            &\Miff \quad \exists w \in \spec{\M} \logdot x \Rb w \text{ and } w \in \pi \nu \chi \\
            &\Miff \quad \exists w \in \spec{\M} \logdot x \Rb w \text{ and } \nu \chi \in w  \\
            &\Miff \quad \exists w \in \spec{\M} \logdot x \Rb w \text{ and } \tup{\spec{\M}, w} \models  \chi  \\
            &\Miff \quad \tup{\spec{\M}, x} \models \bDia \chi
        \end{align*}
    \end{prf}
    Indeed, the algebraic truth lemma implies that an algebraic model has a theory identical to that of its dual relational model.
    \begin{lem}\label{lem:spec-truth-preserving-reflecting}
        Given an algebraic model $\M$ and a formula $\phi \in \Lt$,
            \[ \M \models \phi \quad \Miff \quad \spec{\M} \models \phi. \]
    \end{lem}
    \begin{prf}
        Let $\M := \tup{\A, \nu}$ and observe the following.
        \begin{align*}
            \M \models \phi \quad &\Miff \quad \nu\phi = 1 \\
            &\Miff \quad \pi\nu\phi = \pi 1 = \spec{\M} \quad (\text{because $\pi$ is injective}) \\
            &\Miff \quad \forall x \in \spec{\M} \logdot x \in \pi\nu\phi \\
            &\Miff \quad \forall x \in \spec{\M} \logdot \tup{\spec{\M}, x} \models \phi \quad (\text{by Lemma \ref{lem:truth-lemma}}) \\
            &\Miff \quad \spec{\M} \models \phi
        \end{align*}
    \end{prf}
    \par
    We now embark on a study of filtration on temporal transits.
    This, in conjunction with our truth-preserving/reflecting means of transforming algebraic models into relational models, will allow us to establish the relational and algebraic FMP for $\tHC$.
    \begin{con}[Smallest relational filtration]\label{con:smallest-filtration}
        Given a relational model $\M := \tup{\X, \nu}$ and a finite, subformula-closed set $\Sigma \sub \Lt$, we construct a \textit{finite} relational model $\M_{\Sigma}$.
        ~\\~\\
        We define the binary relation $\sim$ on $\X$ as follows.
            \[ x \sim y \quad \Diff \quad \forall \phi \in \Sigma \logdot \tup{\M, x} \models \phi \Leftrightarrow \tup{\M, y} \models \phi \]
        It is clear that $\sim$ is an equivalence relation.
        We let $X_{\Sigma} := X/{\sim}$ and we use the shorthand $[x] := [x]_{\sim} = \compin{y}{X}{x \sim y}$.
        ~\\~\\
        We then define
        \begin{align*}
            [x] \rbSig [w] \quad &\Diff \quad \exists x', w' \in X \logdot x \sim x' \Rb w' \sim w \\
            [x] \rfSig [y] \quad &\Diff \quad \exists x', y' \in X \logdot x \sim x' \Rf y' \sim y \\
            [x] \precSig [y] \quad &\Diff \quad \exists x', y' \in X \logdot x \sim x' \leq y' \sim y
        \end{align*}
        and define ${\RbSig}$, ${\RfSig}$, and ${\leqSig}$ to be the transitive closures of ${\rbSig}$, ${\rfSig}$, and ${\precSig}$ respectively.
        ~\\~\\
        Next we define $\nuSig: \Sigma \to X_{\Sigma}$ by the rule
            \[ p \mapstoo \compin{[x]}{X_{\Sigma}}{x \in \nu p}. \]
        Finally, we define
            \[ \X_{\Sigma} := \tup{X_{\Sigma}, {\RbSig}, {\RfSig}, {\leqSig}} \Quad \M_{\Sigma} := \tup{\X_{\Sigma}, \nu_{\Sigma}}.\]
        It will be shown in Lemma \ref{lem:smallest-filtration-preserves-tikm} that $\M_{\Sigma}$ is indeed a finite relational model of $\tHC$.
    \end{con}
    \begin{lem}\label{lem:smallest-filtration-preserves-tikm}
        Given a relational model $\M$ and a finite, subformula-closed set $\Sigma \sub \Lt$, $\M_{\Sigma}$ is a \textit{finite} relational model of $\tHC$.
    \end{lem}
    \begin{prf}
        Let $\tup{\X, \nu} := \M$.
        It is well known from \cite[Corollary 5.25]{ChaZak97} that $\M_{\Sigma}$ is finite.
        It is, further, well known that ${\leqSig}$ is posetal on $X_{\Sigma}$ and $\nu_{\Sigma}$ is an intuitionistic valuation on $\X_{\Sigma}$ \cite[\S5.2]{ChaZak97}, so it remains only to check that (1) ${\leqSig}$ is the reflexivisation of $\RfSig$, and (2) ${\RbSig}$ and ${\RfSig}$ are inverses of each other.
        ~\\~\\
        (1)
            Letting $\text{refl}(R)$ and $\text{tran}(R)$ denote the reflexivisation and transitive closure of a relation $R$ respectively, we claim that it suffices to show that $\text{refl}({\rfSig}) = {\precSig}$.
            For if this is the case, then we will have $\text{tran}(\text{refl}({\rfSig})) = \text{tran}({\precSig})$, implying that
                \[ \text{refl}({\RfSig}) = \text{refl}(\text{tran}({\rfSig})) = \text{tran}(\text{refl}({\rfSig})) = \text{tran}({\precSig}) = {\leqSig} \]
            as desired.
            Letting $R:= \text{refl}({\rfSig})$, this amounts to showing that $[x]\R[y] \miff [x]\precSig[y]$.
            ($\Rightarrow$)
                Given $[x]\R[y]$, we distinguish the cases where $[x] = [y]$ and where $[x] \neq [y]$.
                In the former case, we have $[x]\precSig [y]$ as desired due to the reflexivity of ${\precSig}$.
                In the latter case, we have $[x]\rfSig [y]$, implying that $x \sim x' \Rf y' \sim y$, further implying that $x \sim x' \leq y' \sim y$, finally implying that $[x]\precSig[y]$ as desired.
            ($\Leftarrow$)
                Given $[x]\precSig[y]$, we have $x \sim x' \leq y' \sim y$.
                We distinguish cases where $x' = y'$ and where $x' \neq y'$.
                In the former case, we have $x \sim x' = y' \sim y$, implying that $[x] = [y]$, further implying that $[x] \R [y]$ (as $R$ is reflexive by construction).
                In the latter case, we have $x \sim x' < y' \sim y$, implying that $x \sim x' \Rf y' \sim y$, further implying that $[x] \rfSig [y]$, finally implying that $[x]\R[y]$ as desired.
        ~\\~\\
        \vspace{-1pt}
        (2)
            This is clear as
            \begin{align*}
              [x]\RfSig[y] \quad &\Miff \quad  \exists z_1, \dots, z_n \in X \logdot [x] = [z_1] \rfSig \dots \rfSig [z_n] = [y] \\
              &\Miff \quad \exists z_1, \dots, z_n \in X \logdot [y] = [z_n] \rbSig \dots \rbSig [z_1] = [x] \\
              &\Miff \quad [y] \RbSig [x].
            \end{align*}
    \end{prf}
    Having established that the filtration of a relational model through an appropriate set of formulas is a finite relational model, we turn our focus to the preservation of formulas through this process, which requires us to check some well-known filtration conditions.
    \begin{lem}\label{lem:valid-filtration}
        Given a relational model $\M$ and a subformula-closed set $\Sigma \sub \Lt$, the following hold for all $w, x, y \in \M$.
        \begin{enumerate}
            \item $x \Rf y$ implies $[x] \RfSig [y]$
            \item $[x] \RfSig [y]$ implies $\forall \Box\phi \in \Sigma \logdot x \models \Box \phi \mimp y \models \phi$
            \item $x \Rb w$ implies $[x] \RbSig [w]$
            \item $[x] \RbSig [w]$ implies $\forall \bDia \phi \in \Sigma \logdot w \models \phi \mimp x \models \bDia \phi$
            \item $x \leq y$ implies $[x] \leqSig [y]$
            \item $[x] \leqSig [y]$ implies $\forall \phi \in \Sigma \logdot x \models \phi \mimp y \models \phi$
        \end{enumerate}
    \end{lem}
    \begin{prf}
        (1), (3), and (5) follow by definition and (6) is well known from \cite[\S5.2]{ChaZak97}.
        Further, (4) follows from a symmetrical argument to (2), so we prove only (2).
        Given $[x] \RfSig [y]$, we have a finite path through ${\rfSig}$ from $[x]$ to $[y]$ (because the former relation was defined to be the transitive closure of the latter).
        If we let $z_1 := x$ and $z_n := y$, then we have
            \[ [z_1] \rfSig [z_2] \rfSig \dots \rfSig [z_{n-1}] \rfSig [z_n]. \]
        Recalling that
            \[ [z_i] \rfSig [z_{i+1}] \quad :\Longleftrightarrow \quad \exists z_{i}'', z_{i+1}' \in X \logdot z_i \sim z_i'' \Rf z_{i+1}' \sim z_{i+1}, \]
        (the standard smallest filtration), we have the following (where dashed lines represent equivalence up to ${\sim}$).
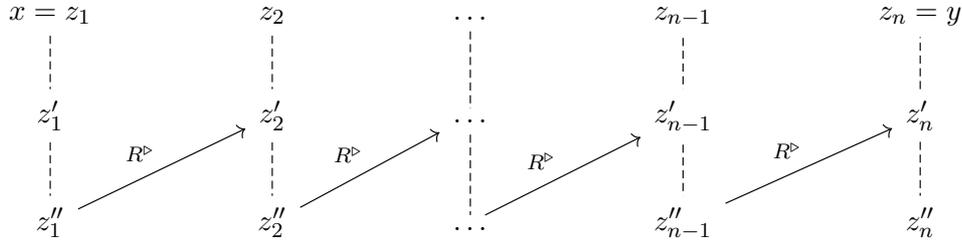
        \[\begin{tikzcd}[column sep=scriptsize,row sep=small]
            {x = z_1} && {z_2} && \dots && {z_{n-1}} && {z_{n} = y} \\
            \\
            {z_1'} && {z_2'} && \dots && {z_{n-1}'} && {z_{n}'} \\
            \\
            {z_1''} && {z_2''} && \dots && {z_{n-1}''} && {z_{n}''}
            \arrow[dashed, no head, from=1-1, to=3-1]
            \arrow[dashed, no head, from=1-3, to=3-3]
            \arrow[dashed, no head, from=1-5, to=3-5]
            \arrow[dashed, no head, from=1-7, to=3-7]
            \arrow[dashed, no head, from=3-1, to=5-1]
            \arrow[dashed, no head, from=3-3, to=5-3]
            \arrow[dashed, no head, from=3-5, to=5-5]
            \arrow[dashed, no head, from=3-7, to=5-7]
            \arrow[dashed, no head, from=3-9, to=1-9]
            \arrow[dashed, no head, from=3-9, to=5-9]
            \arrow["{{\Rf}}", from=5-1, to=3-3]
            \arrow["{{\Rf}}", from=5-3, to=3-5]
            \arrow["{{\Rf}}", from=5-5, to=3-7]
            \arrow["{{\Rf}}", from=5-7, to=3-9]
        \end{tikzcd}\]
        Recalling that ${\Rf} \sub {\leq}$, we have the following.
        
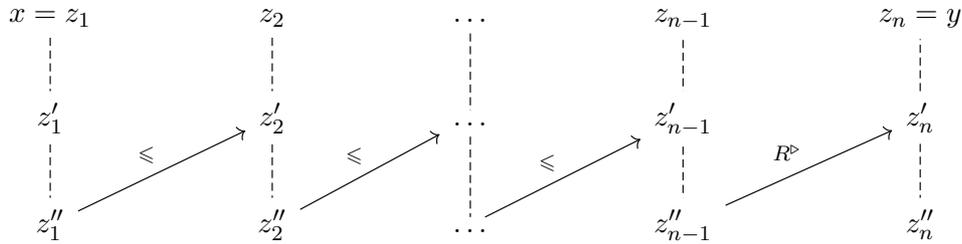
\[\begin{tikzcd}[column sep=scriptsize,row sep=small]
            {x = z_1} && {z_2} && \dots && {z_{n-1}} && {z_{n} = y} \\
            \\
            {z_1'} && {z_2'} && \dots && {z_{n-1}'} && {z_{n}'} \\
            \\
            {z_1''} && {z_2''} && \dots && {z_{n-1}''} && {z_{n}''}
            \arrow[dashed, no head, from=1-1, to=3-1]
            \arrow[dashed, no head, from=1-3, to=3-3]
            \arrow[dashed, no head, from=1-5, to=3-5]
            \arrow[dashed, no head, from=1-7, to=3-7]
            \arrow[dashed, no head, from=3-1, to=5-1]
            \arrow[dashed, no head, from=3-3, to=5-3]
            \arrow[dashed, no head, from=3-5, to=5-5]
            \arrow[dashed, no head, from=3-7, to=5-7]
            \arrow[dashed, no head, from=3-9, to=1-9]
            \arrow[dashed, no head, from=3-9, to=5-9]
            \arrow["{{\leq}}", from=5-1, to=3-3]
            \arrow["{{\leq}}", from=5-3, to=3-5]
            \arrow["{{\leq}}", from=5-5, to=3-7]
            \arrow["{{\Rf}}", from=5-7, to=3-9]
        \end{tikzcd}\]
        So since $x = z_1 \models \Box \phi$, we have $z_1'' \models \Box \phi$ (because $z_1 \sim z_1''$), implying $z_2' \models \Box \phi$ (because $\nu\Box\phi$ is an upset), implying $z_2'' \models \Box \phi$, eventually implying $z_{n-1}'' \models \Box \phi$.
        Since $z_{n-1}'' \Rf z_n$, this implies that $z_{n}' \models \phi$, implying that $y = z_{n} \models \phi$ as desired.
    \end{prf}
    Having addressed the filtration conditions, we can state the \mbox{« filtration} \mbox{lemma »} for the smallest relational filtration.
    \begin{lem}[Relational filtration lemma]\label{lem:filtration}
        Given a relational model $\M$, a subformula-closed set $\Sigma \sub \Lt$, a point $x \in \M$, and $\phi \in \Sigma$,
            \[ \Tup{\M, x} \models \phi \quad \Miff \quad \Tup{\M_{\Sigma}, [x]} \models \phi.  \]
    \end{lem}
    \begin{prf}
        We proceed via induction on the shape of $\phi$.
        Now if $\phi$ has non-modal semantics, this follows either by definition or by a trivial argument, so we consider the cases ${\imp}$, ${\bDia}$, and ${\Box}$.
        Having shown in Lemma \ref{lem:valid-filtration} that ${\RbSig}$ and ${\RfSig}$ are valid \textit{modal filtrations} and that ${\leqSig}$ is a valid \textit{intuitionistic filtration} (both defined in \cite[p.\ 140]{ChaZak97}), this follows from the \mbox{« Filtration} \mbox{Theorem »} of \cite[Theorem 5.23]{ChaZak97}.
    \end{prf}
    \par
    We can now establish the relational and algebraic FMP.
    \begin{thm}[Relational FMP]\label{thm:thc-soco-ttran-fin}
        $\tHC \soco \tTran_{\fin}$
    \end{thm}
    \begin{prf}
        Given $\phi \notin \tHC$, completeness with respect to $\tTran$ implies that there is some relational model $\M := \tup{\X, \nu}$ such that $\M \refutes \phi$, implying that there exists some $x \in \M$ such that $\tup{\M, x} \refutes \phi$.
        If we let
            \[ \Sigma := \compin{\chi}{\Lt}{\chi \text{ is a subformula of } \phi}, \]
        then Lemma \ref{lem:filtration} implies that $\tup{\M_{\Sigma},[x]} \refutes \phi$, implying $\M_{\Sigma} \refutes \phi$, further implying, by Lemma \ref{lem:smallest-filtration-preserves-tikm}, that $\tTran_{\fin} \refutes \phi$.
    \end{prf}
    \begin{thm}[Algebraic FMP]\label{thm:thc-soco-tha-fin}
        $\tHC \soco \tHA_{\fin}$
    \end{thm}
    \begin{prf}
        Given $\phi \notin \tHC$, Theorem \ref{thm:thc-soco-ttran-fin} implies that there is some finite relational model $\M := \tup{\X, \nu}$ such that $\M \refutes \phi$, implying that there exists some $x \in \M$ such that $\tup{\M, x} \refutes \phi$, implying that $\nu\phi \neq \X$.
        Now simply observe that since $\X$ is finite, all upsets are clopen, implying that $\nu : \Prop \to \Up(\X)$ is a well-defined algebraic valuation on $\clop{\X}$ and we have $\nu\phi \neq \X = 1_{\clop{\X}}$.
        (Note that one must confirm that $\nu: \Prop \to \clop{\X}$ extends to a homomorphism $\Term \to \clop{\X}$, but this is quite trivial to check given the comment after Definition \ref{dfn:dual-relational-model} and the fact that all sets are clopen.)
        Since $\clop{\X} \in \tHA_{\fin}$ and $\tup{\clop{\X}, \nu} \refutes \phi$, we can conclude that $\tHA_{\fin} \refutes \phi$.
    \end{prf}
    \par
    The algebraic FMP, in conjunction with a well-known result from universal algebra, allows us to establish an even stronger completeness result with respect to \textit{finite subdirectly-irreducible} algebras.
    \begin{thm}\label{thm:thc-soco-tha-fsi}
        $\tHC \soco \tHA_{\fsi}$.
    \end{thm}
    \begin{prf}
        It is well known that every finite element of a variety $\V$ is isomorphic to the subdirect product of a finite set of finite subdirectly-irreducible algebras in $\V$ \cite[Corollary 3.25]{Ber12}.
        Recalling that $\tHA$ has an equational definition, and is, therefore, a variety, this implies that for all $\A \in \tHA_{\fin}$, we have $\A \iso \B \leq \prod_{i=1}^{n}\C_i$ for some $\bindset{\C_i}{i=1}{n} \sub \tHA_{\fsi}$.
        It is also well known that the relevant class operators $\mathsf{I}$ and $\mathsf{P}_{\mathsf{S}}$ are truth-preserving, implying that if $\tHA \contains \A \refutes \phi$, there is some $\C_{k}$ such that $\C_{k} \refutes \phi$.
        Completeness then follows from the following argument : $\phi \notin \tHC$ implies $\tHA_{\fin} \contains \A \refutes \phi$ (by Theorem \ref{thm:thc-soco-tha-fin}), implying $\B \refutes \phi$, further implying $\tHA_{\fsi} \contains \C_{k} \refutes \phi$ and, therefore, $\tHA_{\fsi} \refutes \phi$.
    \end{prf}
    \par
    In light of Theorem \ref{thm:thc-soco-tha-fsi}, we can apply Corollary \ref{cor:si-characterisation-fin} to prove our final completeness result.
    \begin{thm}\label{thm:thc-soco-fin-z-rooted-ttran}
        $\tHC ~\soco$ The class of finite $Z$-rooted temporal transits.
    \end{thm}
    \begin{prf}
        Let the above-described class be denoted by $\mathbf{Z}$.
        Given $\phi \notin \tHC$, Theorem \ref{thm:thc-soco-tha-fsi} implies that there is some algebraic model $\M := \tup{\A, \nu}$ such that $\A \in \tHA_{\fsi}$ and $\M \refutes \phi$.
        Lemma \ref{lem:spec-well-defined-tha-obj} implies that $\fram{\A}$ is a finite temporal transit and Corollary \ref{cor:si-characterisation-fin} implies that $\fram{\A}$ is $Z$-rooted, so we have $\fram{\A} \in \mathbf{Z}$.
        By Lemma \ref{lem:spec-truth-preserving-reflecting}, we have $\spec{\M} \refutes \phi$.
        Combining these two facts, we can conclude that $\mathbf{Z} \refutes \phi$.
    \end{prf}
\section*{Acknowledgements}
    I would like to thank Rodrigo Nicolau Almeida and Nick Bezhanishvili for their invaluable guidance while supervising my MSc Logic thesis project \cite{Alv24} at the Institute for Logic, Language, and Computation, during which the above results were obtained.
    I would also like to acknowledge Phridon Alshibaia and Guram Bezhanishvili, who laid the groundwork for a duality-theoretic study of the temporal Heyting calculus \cite{Als14}.
    Finally, I would like to thank an anonymous \href{https://link.springer.com/journal/12}{\textit{Algebra Universalis}} reviewer for their thorough and insightful feedback.

\end{document}